\documentclass[11pt]{article}
\usepackage{amsmath,latexsym,epsf,mathabx}
\usepackage{changepage}
\begin{document}

\newcommand{\nc}{\newcommand}
\def\PP#1#2#3{{\mathrm{Pres}}^{#1}_{#2}{#3}\setcounter{equation}{0}}
\def\ns{$n$-star}\setcounter{equation}{0}
\def\nt{$n$-tilting}\setcounter{equation}{0}
\def\Ht#1#2#3{{{\mathrm{Hom}}_{#1}({#2},{#3})}\setcounter{equation}{0}}
\def\qp#1{{${(#1)}$-quasi-projective}\setcounter{equation}{0}}
\def\mr#1{{{\mathrm{#1}}}\setcounter{equation}{0}}
\def\mc#1{{{\mathcal{#1}}}\setcounter{equation}{0}}
\def\HD{\mr{Hom}}
\def\HC{\mr{Hom}_{\mc{C}}}
\def\AdT{\mr{Add}_{\mc{T}}}
\def\adT{\mr{add}_{\mc{T}}}
\def\Kb{\mc{K}^b(\mr{Proj}R)}
\def\kb{\mc{K}^b(\mc{P}_R)}
\def\AdpC{\mr{Adp}_{\mc{C}}}
\def\AdpD{\mr{Adp}_{\mc{D}}}
\newtheorem{Th}{Theorem}[section]
\newtheorem{Def}[Th]{Definition}
\newtheorem{Lem}[Th]{Lemma}
\newtheorem{Pro}[Th]{Proposition}
\newtheorem{Cor}[Th]{Corollary}
\newtheorem{Rem}[Th]{Remark}
\newtheorem{Exm}[Th]{Example}
\newtheorem{Sc}[Th]{}
\def\Pf#1{{\noindent\bf Proof}.\setcounter{equation}{0}}
\def\>#1{{ $\Rightarrow$ }\setcounter{equation}{0}}
\def\<>#1{{ $\Leftrightarrow$ }\setcounter{equation}{0}}
\def\bskip#1{{ \vskip 20pt }\setcounter{equation}{0}}
\def\sskip#1{{ \vskip 5pt }\setcounter{equation}{0}}
\def\mskip#1{{ \vskip 10pt }\setcounter{equation}{0}}
\def\bg#1{\begin{#1}\setcounter{equation}{0}}
\def\ed#1{\end{#1}\setcounter{equation}{0}}
\def\KET{T^{^F\bot}\setcounter{equation}{0}}
\def\KEC{C^{\bot}\setcounter{equation}{0}}

\def\jze{{ \begin{pmatrix} 0 & 0 \\ 1 & 0 \end{pmatrix}}\setcounter{equation}{0}}
\def\hjz#1#2{{ \begin{pmatrix} {#1} & {#2} \end{pmatrix}}\setcounter{equation}{0}}
\def\ljz#1#2{{  \begin{pmatrix} {#1} \\ {#2} \end{pmatrix}}\setcounter{equation}{0}}
\def\jz#1#2#3#4{{  \begin{pmatrix} {#1} & {#2} \\ {#3} & {#4} \end{pmatrix}}\setcounter{equation}{0}}


\title{\bf  Cosilting complexes and AIR-cotilting modules\thanks{Supported by the National Science Foundation of China (Grant No. 11371196) and the National Science Foundation for Distinguished Young Scholars of Jiangsu Province (Grant No. BK2012044) and
a project funded by the Priority Academic Program Development of Jiangsu Higher Education Institutions.}}

\footnotetext{
E-mail:~zhangpeiyu2011@163.com,~weijiaqun@njnu.edu.cn}
\smallskip
\author{\small Peiyu Zhang, Jiaqun Wei\\
\small Institute of Mathematics, School of Mathematics Sciences\\
\small Nanjing Normal University, Nanjing \rm210023 China}
\date{}
\maketitle
\baselineskip 15pt
%
%
\begin{abstract}
\vskip 10pt%
%
%
We introduce and study the new concepts of cosilting complexes, cosilting modules and AIR-cotilting modules. We prove that the three concepts AIR-cotilting modules, cosilting modules and quasi-cotilting modules coincide with each other, in contrast with the dual fact that AIR-tilting modules, silting modules and quasi-tilting modules are different. Further, we show that there are bijections between the following four classes (1) equivalent classes of AIR-cotilting (resp., cosilting, quasi-cotilting) modules, (2) equivalent classes of 2-term cosilting complexes, $(3)$ torsion-free cover classes and $(4)$ torsion-free special precover classes. We also extend a classical result of Auslander and Reiten on the correspondence between certain contravariantly finite subcategories and cotilting modules to the case of cosilting complexes.

\mskip\





\end{abstract}
%
\vskip 30pt

\section{Introduction}

%
%
%
\hskip 18pt
The tilting theory is well known, and plays an important role in the representation theory of
Artin algebra.
The classical notion of tilting and cotilting modules was first considered in the case
of finite dimensional algebras by Brenner and Butler \cite{SBMB} and by Happel and Ringel \cite{HR}.
Cotilting theory (for arbitrary modules over arbitrary unital rings) extends Morita duality in analogy to the way tilting theory extends Morita equivalence. In particular, cotilting modules generalize injective cogenerators similarly as tilting modules generalize progenerators.
Later, many scholars have done a lot of research on the tilting theory and cotilting theory, for instance
\cite{AC, ATT, AR, Bazz-nt, CDT, RCJT, YMIY, Wnt, Wns} and so on.

The silting theory seems to be the tilting theory in the level of derived categories (while the tilting complexes play the role of progenerators).  Silting complexes were first introduced by Keller and Vossieck \cite{BKDV} to study t-structures in the bounded derived category of representations of Dynkin quivers. Beginning with \cite{ATIO}, such objects were recently shown to have various nice properties \cite{{KY}, {MSSS}}. The results in \cite{Wst} show that silting complexes (i.e., semi-tilting complexes in \cite{Wst}) have similar properties as that tilting modules have in the module categories. The recent paper by Buan and Zhou \cite{BZ} also shows that it is reasonable to see the silting theory as the tilting theory in the level of derived categories.

The $\tau$-tilting theory recently introduced by Adachi, Iyama and Reiten \cite{AIR} is an important generalization of the classical tilting theory. In particular, it was shown that support $\tau$-tilting modules have close relations with 2-term silting complexes and cluster-tilting objects \cite{AIR}. In \cite{AMV}, the authors introduce silting modules as a generalization of support $\tau$-tilting modules over arbitrary rings and modules. We note that there is also another little different generalization of support $\tau$-tilting modules over  arbitrary rings, called large support $\tau$-tilting modules, which was introduced by the second author \cite{Wlt}.

\vskip 10pt

In this paper, we concentre on the dual case, i.e., the correspondent cotilting parts to the above achievements. In the level of finitely generated modules over artin algebras, such dual cases proceed very well. So we only consider the dual over arbitrary rings and modules. As one will see, there are many interesting properties in such case.

Let us briefly introduce the contents and main results of this paper in the following.

After the introduction in Section 1, Section 2 is devoted to studying cosilting complexes. Namely,  a complex $T$ over a ring $R$ is cosilting if it satisfies the following three conditions:

\hangafter 1 \hangindent 35pt
(1) $T\in \mc{K}^b(\mr{Inj}R)$,

\hangafter 1 \hangindent 35pt
 (2) $T$ is prod-semi-selforthogonal and,

\hangafter 1 \hangindent 35pt
 (3) $\mc{K}^b(\mr{Inj}R)$ is just the smallest triangulated subcategory containing $\langle \AdpD{T}\rangle$, where $\AdpD{T}$ denotes the class of complexes isomorphic in the derived category $\mc{D}(\mathrm{Mod}R)$ to a direct summand of some direct products of copies of $T$.

 It is clear that a cotilting complex  \cite{Buan} is cosilting. We show that an $R$-module is a cotilting module 
  if and only if it is isomorphic in the derived category to a cosilting complex. Some characterizations of cosilting complexes are obtained. In particular, we extend a simple characterization of cotilting modules  \cite{Bazz-nt} to cosilting complexes  (Theorem \ref{c2}). In \cite{AR},
Auslander and Reiten showed that, over an artin algebra, there is a one-one correspondence between certain contravariantly finite subcategories and basic cotilting modules. The result was extended to the derived category of artin algebras by Buan \cite{Buan}, where the author proved that there is a one-one correspondence between basic cotilting complexes and certain contravariantly finite subcategories of the derived category. Here, we further extend the result to cosilting complexes and to arbitrary rings (Theorem \ref{TIAO}).

In Section 3, we study quasi-cotilting modules and cosilting modules. In the tilting case, it is known that a silting module is always a finendo quasi-tilting module but the converse is not true in general \cite{Trf}. However, in the dual case, we see that quasi-cotilting modules are always cofinendo and that they are also pure-injective \cite{ZW}. We show that cosilting modules are always quasi-cotilting modules, and consequently, cosilting modules are pure injective and cofinendo. Interesting properties and characterizations of cosilting modules are also given in this section.

In Section 4, we introduce and study AIR-cotilting modules. We call an $R$-module $M$ AIR-tilting if it is large support $\tau$-tilting in sense of \cite{Wlt}, i.e., it satisfies the following two conditions:


\hangafter 1 \hangindent 35pt
 (1) there is an exact sequence $P_1\to^{f} P_0\to M\to 0$ with $P_1, P_0$ projective such that $\mathrm{Hom}(f,M^{(X)})$ is surjective for any set $X$ and,

\hangafter 1 \hangindent 35pt
  (2) there is an exact sequence $R\to^g M_0\to M_1\to 0$ with $M_0, M_1\in \mathrm{Add}M$ such that $\mathrm{Hom}(g,M^{(X)})$ is surjective for any set $X$.


\noindent  Dually, we call an $R$-module $M$ AIR-cotilting if  it satisfies the following two conditions:

\hangafter 1 \hangindent 35pt
(1) there is an exact sequence $0\to T\to I_0\to^{f} I_1$ with $I_0, I_1$ injective such that $\mathrm{Hom}(T^{X},f)$ is surjective for any set $X$ and,

\hangafter 1 \hangindent 35pt
 (2) there is an exact sequence $0\to T_1\to T_0\to^{g} Q$ with $T_0, T_1\in \mathrm{Adp}T$ and $Q$ an injective cogenerator such that $\mathrm{Hom}(T^{X},g)$ is surjective for any set $X$.

 Clearly, 1-tilting modules are AIR-tilting modules and 1-cotilting modules are AIR-cotilting modules. In the tilting case, a silting module is always AIR-tilting and an AIR-tilting module can be completed to a silting module \cite{AMV, Wlt}. But it is a question if these two notions are the same in general (they are the same in the scope of finitely generated modules over artin algebras). It is also known that AIR-tilting modules are finendo quasi-tilting. But the converse is not true in general \cite{Trf}. However, in the dual case, we prove  that AIR-cotilting modules coincide with cosilting modules, as well as quasi-cotilting modules (Theorem \ref{31}). Moreover, we also show that there is a 1-1 correspondence between equivalent classes of AIR-cotilting modules and 2-term cosilting complexes (Theorem \ref{AIRcs}).

Summarized, we obtain the following main results.

\bg{Th}\label{}%
There is a one-one correspondence, given by $u: T \mapsto$ ${^{\bot_{i>0}}T}$,
between equivalent classes of cosilting complexes in $\mathcal{D}^{\geq}$
and subcategories $\mc{T}\subseteq \mathcal{D}^{\geq}$ which is specially
contravariantly finite in $\mathcal{D}^{+}$, resolving and closed under
products such that $\widehat{\mathcal{T}}=\mathcal{D}^{+}$.
\ed{Th}%

\bg{Th}\label{}%

Let $R$ be a ring and $M$ be an $R$-module. The following statements are equivalent:

$(1)$ $M$ is AIR-cotilting;

$(2)$ $M$ is quasi-cotilting.

$(3)$ $M$ is cosilting.
\ed{Th}%

\bg{Th}\label{}%

There are bijections between

$(1)$ equivalent classes of AIR-cotilting (resp., cosilting, quasi-cotilting) modules;

$(2)$ equivalent classes of 2-term cosilting complexes;

$(3)$ torsion-free cover classes and,

$(4)$ torsion-free special precover classes.

\ed{Th}
\mskip\

Throughout this paper, $R$ will denote an associative ring with identity and we mainly work on the category of left $R$-modules which is denoted by $\mr{Mod}R$.
We denote by $\mathrm{Inj}R$ (resp., $\mathrm{Proj}R$) the class of all injective (resp., projective) $R$-modules.
The notations $\mc{K}^b(\mr{Inj}R)$ (resp., $\mc{K}^b(\mr{Proj}R)$) denotes the homotopy category of bounded (always cochain) complexes of
injective (resp., projective) modules.
The unbounded derived category of $\mathrm{Mod}R$ will be denoted by $\mathcal{D}(R)$, or simply $\mathcal{D}$, with [1] shift functor. We denoted by $\mathcal{D}^{\geq}$ the subcategory of complexes whose homologies
are concentrated on non-negative terms. We use $\mc{D}^{+}$ to denote the subcategory of $\mc{D}$ consists of bounded-below complexes.

Note that $\mc{D}$ is a triangulated category and $\mc{K}^b(\mr{Inj}R)$, $\mc{K}^b(\mr{Proj}R)$, $\mc{D}^{+}$ are all full triangulated subcategories of $\mc{D}$. We refer to Happel's paper \cite{Hp-b} for more details on derived categories and triangulated categories.


\vskip 30pt
\section {Subcategories of the derived category and cosilting complexes}
%
%

%
\

In this section we study the dual of silting complexes and give various characterizations of cosilting complexes. In particular, we extend a result of Bazzoni \cite{Bazz-nt} and establish a one-one correspondence between certain subcategories of the derived category $\mc{D}$ and the equivalent classes of  cotilting complexes stemming from Auslander and Reiten \cite{AR}.

We begin with some basic notations and some useful facts in general triangulated categories.

Let $\mc{C}$ be a triangulated category with [1] the shift functor.
Assume that $\mc{B}$ is a full subcategory of $\mc{C}$. Recall that $\mc{B}$ is closed under extension if for any triangle $X\to Y\to Z\to$ in  $\mc{C}$ with
$X, Z \in \mc{C}$, we have $V\in \mc{C}$. The subcategory
$\mc{B}$ is resolving (resp., coresolving) if it is
closed under extension and under the functor $[-1]$ (resp., $[1]$). It is
easy to prove that $\mc{B}$ is resolving (resp., coresolving) if and
only if for any triangle $X\to Y\to Z\to$ (resp.,  $Z\to Y\to
X\to$) in $\mc{B}$ with $Z\in\mc{B}$, one has that `\
$X\in\mc{B}\Leftrightarrow Y\in\mc{B}$\ '.

We say that an object $M\in \mc{C}$ has a $\mc{B}$-resolution (resp., $\mc{B}$-coresolution) with the length
at most $m$ ($m\ge 0$), if there are
triangles $M_{i+1}\to X_i\to M_i\to$ (resp., $M_i\to X_i\to
M_{i+1}\to$) with $0\le i\le m$ such that $M_0=M$, $M_{m+1}=0$
and each $X_i\in \mc{B}$. In the case, we denoted by $\mc{B}$-res.dim$(L)\le m$
(resp., $\mc{B}$-cores.dim$(L)\le m$). One may compare such notions with the usual finite resolutions and coresolutions respectively in the module category.

Associated with a subcategory $\mc{B}$, we have the following notations which are widely used in the tilting theory (see for instance \cite{AR}),
where $n\ge 0$ and $m$ is an integer.

\sskip\

\hskip 25pt $(\widehat{\mc{B}})_n=\{L\in\mc{C}\ |\
\mc{B}$-res.dim$(L)\le n\}$.

\hskip 25pt $(\widecheck{\mc{B}})_n=\{L\in\mc{C}\ |\
\mc{B}$-cores.dim$(L)\le n\}$.

\hskip 25pt $\widehat{\mc{B}}\ \ \ \ =\{L\in\mc{C}\ |\ L\in
(\widehat{\mc{B}})_n$ for some $n\}$.

\hskip 25pt $\widecheck{\mc{B}}\ \ \ \ =\{L\in\mc{C}\ |\ L\in
(\widecheck{\mc{B}})_n$ for some $n\}$.

\hskip 25pt ${\mc{B}}^{\bot_{i\neq 0}}\ =\{N\in \mc{C}\ |\
\HD(M,N[i])=0$ for all $M\in \mc{B}$ and all $i\neq 0\}$.

\hskip 25pt $^{\bot_{i\neq 0}}{\mc{B}}\ =\{N\in \mc{C}\ |\
\HD(N,M[i])=0$ for all $M\in \mc{B}$ and all $i\neq 0\}$.

\hskip 25pt ${\mc{B}}^{\bot_{i>m}}=\{N\in \mc{C}\ |\
\HD(M,N[i])=0$ for all $M\in \mc{B}$ and all $i>m\}$.

\hskip 25pt $^{\bot_{i>m}}{\mc{B}}=\{N\in \mc{C}\ |\
\HD(N,M[i])=0$ for all $M\in \mc{B}$ and all $i>m\}$.

\hskip 25pt ${\mc{B}}^{\bot_{i\gg 0}}\ =\{N\in \mc{C}\ |\
N\in{\mc{B}}^{\bot_{i>m}}$ for some $m\}$.

\sskip\ Note that ${\mc{B}}^{\bot_{i>m}}$ (resp.,
$^{\bot_{i>m}}\mc{B}$) is coresolving (resp., resolving) and
closed under direct summands and that ${\mc{B}}^{\bot_{i\gg 0}}$  is a triangulated subcategory of $\mc{C}$.

The subcategory $\mc{B}$ is said to be semi-selforthogonal (resp., selforthogonal) if $\mc{B}\subseteq {\mc{B}}^{\bot_{i>m}}$ (resp., $\mc{B}\subseteq {\mc{B}}^{\bot_{i\neq 0}}$). For instance, both subcategories $\mr{Proj}R$ and $\mr{Inj}R$ are selforthogonal.

\vskip 10pt

{\it In the following results of this section, we always assume that $\mc{B}$ is semi-selforthogonal and that $\mc{B}$ is additively closed} (i.e., $\mc{B}=\mathrm{add}_{\mc{C}}\mc{B}$ where $\mathrm{add}_{\mc{C}}\mc{B}$ denotes the subcategory of all objects in $\mc{C}$ which are isomorphic to a direct summand of finite direct sums of copies of objects in $\mc{B}$).

Associated with the subcategory $\mc{B}$, we also have the following two useful subcategories  which are again widely used in the tilting theory (see for instance \cite{AR}).

\bg{verse}
 ${\mc{X}_{\mc{B}}}\ \ =\{N\in {^{\bot_{i>0}}}\mc{B}\ |$ there are
 triangles $N_{i}\to B_i\to N_{i+1}\to$  such that
$N_0=N$, $N_i\in {^{\bot_{i>0}}\mc{B}}$ and $B_i\in\mc{B}$ for
all $i\ge 0\}$.

${_{\mc{B}}\mc{X}}\ =\{N\in {\mc{B}}^{\bot_{i>0}}\ |$ there are
 triangles $N_{i+1}\to B_i\to N_i\to$  such that $N_0=N$,
$N_i\in {\mc{B}}^{\bot_{i>0}}$ and $B_i\in\mc{B}$ for all $i\ge 0\}$.

\ed{verse}

We summarize some results on subcategories associated with $\mc{B}$ in the following, where $\langle\mc{B}\rangle$ denotes the smallest triangulated subcategory containing $\mc{B}$. We refer to \cite{Wst} for their proofs.
%
%
%

\bg{Pro}\label{TIRA}\label{IRAC}\label{AACU}\label{OIIE}%
Let $\mc{B}$ be a semi-selforthogonal subcategory of a triangulated category $\mc{C}$ such that $\mc{B}$ is additively closed. Then

$(1)$ The three subcategories $\widecheck{\mc{B}}\subseteq {\mc{X}_{\mc{B}}}\subseteq {^{\bot_{i>0}}\mc{B}}$  is resolving and closed under direct summands.

$(2)$ The three subcategories $\widehat{\mc{B}}\subseteq {_{\mc{B}}\mc{X}}\subseteq {\mc{B}^{\bot_{i>0}}}$  is coresolving and closed under direct summands.

$(3)$ $\mc{B}=\widecheck{\mc{B}}\bigcap {\mc{B}}^{\bot_{i>0}}$=$\widehat{\mc{B}}\bigcap {^{\bot_{i>0}}\mc{B}}$.

\hangafter 1 \hangindent 35pt
$(4)$ $(\widecheck{\mc{B}})_n={\mc{X}_{\mc{B}}}\bigcap
({\mc{X}_{\mc{B}}})^{\bot_{i>n}}={\mc{X}_{\mc{B}}}\bigcap
({^{\bot_{i>0}}\mc{B}})^{\bot_{i>n}}$. In particular, it is closed under extensions and direct summands.

 $(5)$ $(\widehat{\mc{B}})_n={_{\mc{B}}\mc{X}}\bigcap
{^{\bot_{i>n}}({_{\mc{B}}\mc{X}})}$
$={_{\mc{B}}\mc{X}}\bigcap
{^{\bot_{i>n}}({\mc{B}}^{\bot_{i>0}})}$. In particular, it is closed under extensions and direct summands.

$(6)$ The following three subcategories coincide with each other.

\begin{adjustwidth}{30pt}{0cm}

\noindent $(i)$ $\langle\mc{B}\rangle$: the smallest triangulated
subcategory containing $\mc{B}$;

\noindent $(ii)$ $(\widehat{\mc{B}})_{-}=\{X\in \mc{C} |$ there exists some
$Y\in \widehat{\mc{B}}$ and some $i\le 0$ such that $X=Y[i]\}$;

\noindent $(iii)$  $(\widecheck{\mc{B}})_{+}=\{X\in \mc{C} |$ there exists some
$Y\in \widecheck{\mc{B}}$ and some $i\ge 0$ such that $X=Y[i]\}$.

\end{adjustwidth}

$(7)$ $\widehat{\mc{B}}={\mc{B}}^{\bot_{i>0}}\bigcap\
\langle\mc{B}\rangle$.

$(8)$ $\widecheck{\mc{B}}={^{\bot_{i>0}}\mc{B}}\bigcap\
\langle\mc{B}\rangle$.

\ed{Pro}

We also need the following results.

\bg{Lem}\label{STAT}%
Suppose that $n\ge 1$ and there are triangles $L_{i}\to M_i\to
L_{i+1}\to$ with each $M_i\in {\mc{X}_{\mc{B}}}$, where $0\le i\le
n-1$. Then there exist $X_n, Y_n\in \mc{C}$ such that

$(1)$ $Y_n\in {\mc{X}_{\mc{B}}}$,

$(2)$ there is a triangle $L_{n}\to X_{n}\to Y_{n}\to$, and

\hangafter 1 \hangindent 35pt
$(3)$  there are triangles $X_{i-1}\to B_{i-1}\to X_{i}\to$ with
each $B_{i-1}\in \mc{B}$, for all $1\le i\le n$, where $X_0=L_0$.
\ed{Lem}

\Pf. We use induction on $n$ to prove this conclusion.

\ \hskip 17pt For $n=1$, there is a triangle $M_{0}\to B_0\to Y_{1}\to$
with $B_0\in \mc{B}$ and $Y_{1}\in {\mc{X}_{\mc{B}}}$ since $M_0\in {\mc{X}_{\mc{B}}}$.
Then we can get the following triangle commutative diagram:
\mskip\

 \setlength{\unitlength}{0.09in}
 \begin{picture}(50,18)

                 \put(18,3.4){\vector(0,-1){2}}
                 \put(27,3.4){\vector(0,-1){2}}
                 \put(35,3.4){\vector(0,-1){2}}

 \put(18,5){\makebox(0,0)[c]{$0$}}
                             \put(21,5){\vector(1,0){2}}
 \put(27,5){\makebox(0,0)[c]{$Y_1$}}
                             \put(30,5){\vector(1,0){2}}
 \put(35,5){\makebox(0,0)[c]{$Y_1$}}
                             \put(37,5){\vector(1,0){2}}

                 \put(18,9){\vector(0,-1){2}}
                 \put(27,9){\vector(0,-1){2}}
                 \put(35,9){\vector(0,-1){2}}

 \put(18,11){\makebox(0,0)[c]{$L_0$}}
                             \put(21,11){\vector(1,0){2}}
 \put(27,11){\makebox(0,0)[c]{$B_{0}$}}
                             \put(30,11){\vector(1,0){2}}
 \put(35,11){\makebox(0,0)[c]{$X_1$}}
                             \put(37,11){\vector(1,0){2}}

                 \put(18,14.5){\vector(0,-1){2}}
                 \put(27,14.5){\vector(0,-1){2}}
                 \put(35,14.5){\vector(0,-1){2}}

 \put(18,16){\makebox(0,0)[c]{$L_0$}}
                              \put(21,16){\vector(1,0){2}}
 \put(27,16){\makebox(0,0)[c]{$M_0$}}
                              \put(30,16){\vector(1,0){2}}
 \put(35,16){\makebox(0,0)[c]{$L_{1}$}}
                              \put(37,16){\vector(1,0){2}}

\end{picture}

Obviously, $X_{1}$ and $Y_{1}$ in diagram above are just the objects we look for.

We suppose that the result holds for $n-1$. Next, we will verify that the result holds for $n$.
According to the known condition, we have the triangle $L_{n-1}\to M_{n-1}\to L_{n}\to$ with $M_{n-1}\in{\mc{X}_{\mc{B}}}$.
Using the induction on $L_{n-1}$, one can obtain some triangles $L_{n-1}\to X_{n-1}\to Y_{n-1}\to$
with $Y_{n-1}\in  {\mc{X}_{\mc{B}}}$ and $X_{i-1}\to B_{i-1}\to X_{i}\to$ with
$B_i\in \mc{B}$, for all $1\le i\le n-1$, where $X_0=L_0$.
Hence we have the following triangle commutative diagram:

 \setlength{\unitlength}{0.09in}
 \begin{picture}(50,20)

                 \put(18,3.4){\vector(0,-1){2}}
                 \put(27,3.4){\vector(0,-1){2}}
                 \put(35,3.4){\vector(0,-1){2}}

 \put(18,5){\makebox(0,0)[c]{$Y_{n-1}$}}
                             \put(21,5){\vector(1,0){2}}
 \put(27,5){\makebox(0,0)[c]{$Y_{n-1}$}}
                             \put(30,5){\vector(1,0){2}}
 \put(35,5){\makebox(0,0)[c]{$0$}}
                             \put(37,5){\vector(1,0){2}}

                 \put(18,9){\vector(0,-1){2}}
                 \put(27,9){\vector(0,-1){2}}
                 \put(35,9){\vector(0,-1){2}}

 \put(18,11){\makebox(0,0)[c]{$X_{n-1}$}}
                             \put(21,11){\vector(1,0){2}}
 \put(27,11){\makebox(0,0)[c]{$H$}}
                             \put(30,11){\vector(1,0){2}}
 \put(35,11){\makebox(0,0)[c]{$L_{n}$}}
                             \put(37,11){\vector(1,0){2}}

                 \put(18,14.5){\vector(0,-1){2}}
                 \put(27,14.5){\vector(0,-1){2}}
                 \put(35,14.5){\vector(0,-1){2}}

 \put(18,16){\makebox(0,0)[c]{$L_{n-1}$}}
                              \put(21,16){\vector(1,0){2}}

 \put(27,16){\makebox(0,0)[c]{$M_{n-1}$}}
                              \put(30,16){\vector(1,0){2}}
 \put(35,16){\makebox(0,0)[c]{$L_{n}$}}
                              \put(37,16){\vector(1,0){2}}

\end{picture}

From the second column in diagram above, we can obtain that $H\in{\mc{X}_{\mc{B}}}$
by Proposition \ref{IRAC}. So one have a triangle $H\to B_{n-1}\to X_{n}\to$ with
$B_{n-1}\in \mc{B}$, $X_{n}\in{\mc{X}_{\mc{B}}}$.
Consequently, one have the following triangle commutative diagram:

 \setlength{\unitlength}{0.09in}
 \begin{picture}(50,20)

                 \put(18,3.4){\vector(0,-1){2}}
                 \put(27,3.4){\vector(0,-1){2}}
                 \put(35,3.4){\vector(0,-1){2}}

 \put(18,5){\makebox(0,0)[c]{$0$}}
                             \put(21,5){\vector(1,0){2}}
 \put(27,5){\makebox(0,0)[c]{$Y_{n}$}}
                             \put(30,5){\vector(1,0){2}}
 \put(35,5){\makebox(0,0)[c]{$Y_{n}$}}
                             \put(37,5){\vector(1,0){2}}

                 \put(18,9){\vector(0,-1){2}}
                 \put(27,9){\vector(0,-1){2}}
                 \put(35,9){\vector(0,-1){2}}

 \put(18,11){\makebox(0,0)[c]{$X_{n-1}$}}
                             \put(21,11){\vector(1,0){2}}
 \put(27,11){\makebox(0,0)[c]{$B_{n-1}$}}
                             \put(30,11){\vector(1,0){2}}
 \put(35,11){\makebox(0,0)[c]{$X_{n}$}}
                             \put(37,11){\vector(1,0){2}}

                 \put(18,14.5){\vector(0,-1){2}}
                 \put(27,14.5){\vector(0,-1){2}}
                 \put(35,14.5){\vector(0,-1){2}}

 \put(18,16){\makebox(0,0)[c]{$X_{n-1}$}}
                              \put(21,16){\vector(1,0){2}}

 \put(27,16){\makebox(0,0)[c]{$H$}}
                              \put(30,16){\vector(1,0){2}}
 \put(35,16){\makebox(0,0)[c]{$L_{n}$}}
                              \put(37,16){\vector(1,0){2}}

\end{picture}

It is easy to see that
$X_{n}$ and $Y_{n}$ from diagram above are just the objects we want.
{\ \hfill$\Box$}

\bg{Cor}\label{FALT}%
For any $L\in (\widehat{{\mc{X}_{\mc{B}}}})_n$, then there are two
triangles $L\to X\to Y\to$ with $X \in(\widehat{\mc{B}})_{n}$,
$Y\in {\mc{X}_{\mc{B}}}$ and $U\to V\to L\to$ with $U \in(\widehat{\mc{B}})_{n-1}$
and $V\in$ $^{\bot_{i>0}}\mc{B}$.
\ed{Cor}%

\Pf. For any $L\in (\widehat{{\mc{X}_{\mc{B}}}})_n$, we have triangles
$L_{i}\to M_{i}\to L_{i+1}\to$ with $M_{i}\in {\mc{X}_{\mc{B}}}$, where $0\le i\le n$,
$L_{0}=0, L_{n+1}=L$. By Lemma \ref{STAT}, we obtain the triangle
$L\to X\to Y\to$ with $X \in(\widehat{\mc{B}})_{n}$,
$Y\in {\mc{X}_{\mc{B}}}$.
Since $X \in(\widehat{\mc{B}})_{n}$, one can get a triangle $U\to B_{0}\to X\to$
with $T_{0}\in \mc{B}$ and $U\in(\widehat{\mc{B}})_{n-1}$.
Then we have the following triangles commutative diagram:

 \setlength{\unitlength}{0.09in}
 \begin{picture}(50,20)

                 \put(18,3.4){\vector(0,-1){2}}
                 \put(27,3.4){\vector(0,-1){2}}
                 \put(35,3.4){\vector(0,-1){2}}

 \put(18,5){\makebox(0,0)[c]{$L$}}
                             \put(21,5){\vector(1,0){2}}
 \put(27,5){\makebox(0,0)[c]{$X$}}
                             \put(30,5){\vector(1,0){2}}
 \put(35,5){\makebox(0,0)[c]{$Y$}}
                             \put(37,5){\vector(1,0){2}}

                 \put(18,9){\vector(0,-1){2}}
                 \put(27,9){\vector(0,-1){2}}
                 \put(35,9){\vector(0,-1){2}}

 \put(18,11){\makebox(0,0)[c]{$V$}}
                             \put(21,11){\vector(1,0){2}}
 \put(27,11){\makebox(0,0)[c]{$B_{0}$}}
                             \put(30,11){\vector(1,0){2}}
 \put(35,11){\makebox(0,0)[c]{$Y$}}
                             \put(37,11){\vector(1,0){2}}

                 \put(18,14.5){\vector(0,-1){2}}
                 \put(27,14.5){\vector(0,-1){2}}
                 \put(35,14.5){\vector(0,-1){2}}

 \put(18,16){\makebox(0,0)[c]{$U$}}
                              \put(21,16){\vector(1,0){2}}

 \put(27,16){\makebox(0,0)[c]{$U$}}
                              \put(30,16){\vector(1,0){2}}
 \put(35,16){\makebox(0,0)[c]{$0$}}
                              \put(37,16){\vector(1,0){2}}

\end{picture}

From the second row in diagram above, one can easy see that $V\in$ $^{\bot_{i>0}}\mc{B}$
since $Y, T_{0}\in$ $^{\bot_{i>0}}\mc{B}$. Hence the triangle
$U\to V\to L\to$ is just what we want.
\ \hfill $\Box$


%
%
%
%
%
%
%
%
%
%

\def\HD{\mr{Hom}_{\mc{C}}}



\def\HD{\mr{Hom}_{\mc{D}}}
\def\Ad{\mr{Add}_{\mc{D}}}
\def\ad{\mr{add}_{\mc{D}}}
\def\D<{\mc{D}^{\le 0}}
\def\Dn{\mc{D}^{-}}

%

\def\KL#1{{^{\bot_{i>0}}(#1)}}
\def\KR#1{{{#1}^{\bot_{i>0}}}}
\def\Kg#1{{{#1}^{\bot_{i\gg 0}}}}
\def\PD{\mr{Pres}_{\D<}^n(\Ad{T})}


\hskip 20pt

Now let $R$ be a ring and $T$ be a complex. Recall that $\AdpD{T}$ denotes the class of complexes isomorphic in the derived category $\mc{D}$ to a direct summand of some direct products of $T$. We say that $T$ is \textit{prod-semi-selforthogonal} if $\AdpD{T}$ is semi-selforthogonal. It is easy to see that $\AdpD{T}$ is additively closed in this case. So the results above applies when we set $\mc{B}=\AdpD{T}$.

We introduce the following definition.


\bg{Def}\label{ACTI}%
 A complex $T$ is said to cosilting if it satisfies the following conditions:

 $(1)$ $T\in \mc{K}^b(\mr{Inj}R)$,

 $(2)$ $T$ is prod-semi-selforthogonal, and

 \hangafter 1 \hangindent 35pt
$(3)$ $\mc{K}^b(\mr{Inj}R) = \langle \AdpD{T}\rangle$, i.e., $\mc{K}^b(\mr{Inj}R)$ coincides with the  smallest triangulated
subcategory containing $\AdpD{T}$.
\ed{Def}%

Now let $Q$ be an injective cogenerator for $\mathrm{Mod}R$. Recall that $\mc{D}^{+}$ is the triangulated subcategory of the derived category $\mc{D}$ consists of bounded-below complexes,
$\mc{K}^b(\mr{Inj}R)$ is homotopy category of bounded complexes of
injective modules. Also, recall that $\mathcal{D}^{\geq}$ is the subcategory of the complexes whose homologies
are concentrated on non-negative terms. It is not difficult to verify that $\mc{D}^{\geq}=$ $^{\bot_{i>0}}Q$
and $\mc{D}^{+}=$ $^{\bot_{i\gg0}}Q$.

%
%
%
%

%


The following result gives a characterization of cosilting complexes.

\bg{Th}\label{ATTA}%
Let $T$ be a complex and Q be an injective cogenerator of $\mathrm{Mod}R$. Up to shifts, we may assume that $T\in \mc{D}^{\geq}$.  Then $T$ is cosilting if
and only if it satisfies the following three conditions:

$(i)$ $T\in \widecheck{\AdpD{Q}}$,

$(ii)$ $T$ is prod-semi-selforthogonal, and

$(iii)$ $Q\in \widehat{\AdpD{T}}$.

\ed{Th}%

\Pf. $\Leftarrow$ Since $T\in \widecheck{\AdpD{Q}}$, there are  triangles $T_{i}\to Q_{i}\to T_{i+1}\to$
with $Q_{i} \in \AdpD{Q}$ for all $0\leq i\leq n$, where $T_{0}=T, T_{n+1}=0$.
It is easy to see that each $T_{i}\in \mc{K}^b(\mr{Inj}R)$, for $i=n, n-1,\cdots\ , 0$, since $\AdpD{Q}=$Inj $R$ $\subseteq \mc{K}^b(\mr{Inj}R)$. In particular, $T\in \mc{K}^b(\mr{Inj}R)$.
Now we need only prove that $\mc{K}^b(\mr{Inj}R) = \langle \AdpD{T}\rangle$.
Note that $\AdpD{Q}\subseteq \widehat{\AdpD{T}}$ by Proposition \ref{AACU}, since $Q\in \widehat{\AdpD{T}}$.
and that $T\in \mc{K}^b(\mr{Inj}R)$, so we have
$$\mc{K}^b(\mr{Inj}R)= \langle \AdpD{Q}\rangle\subseteq\langle \widehat{\AdpD{T}}\rangle
=\langle \AdpD{T} \rangle\subseteq \mc{K}^b(\mr{Inj}R).$$
Thus, $T$ is cosilting.

$\Rightarrow$ Since $T\in \mc{D}^{\geq}=$ $^{\bot_{i>0}}Q$ and $T\in \langle\AdpD{T}\rangle=\mc{K}^b(\mr{Inj}R)
=\langle \AdpD{Q}\rangle$, one can get that $T\in \widecheck{\AdpD{Q}}$ by  Proposition \ref{AACU}.
Note that $\AdpD{T}\subseteq$ $^{\bot_{i>0}}Q$, since $T\in \mc{D}^{\geq}=$  $^{\bot_{i>0}}Q$, so
$Q\in (\AdpD{T})^{\bot_{i>0}}$. Combining with the fact that $Q\in \mc{K}^b(\mr{Inj}R)=\langle \AdpD{T}\rangle$,
we have that $Q\in \widehat{\AdpD{T}}$ by  Proposition \ref{AACU}.
\ \hfill $\Box$

\vskip 10pt

Recall that an $R$-module $T$ is ($n$-)cotilting (see for instance \cite{Bazz-nt}) if it satisfies the following three conditions (1) $\mr{id}T\le n$, i.e., the injective dimension of $T$ is finite, (2) $\mr{Ext}_R^i(T^X,T)=0$ for any $X$ and, (3) there is an exact sequence $0\to T_n\to\cdots\to T_0\to Q\to 0$, where $T_i\in\mr{Adp}_R{T}$ and $Q$ is an injective cogenerator of $\mr{Mod}R$.

\bg{Pro}\label{ATTI}%
Assume that T is an R-module. Then $T$ is a cotilting module if and only if
$T$ is  isomorphic in the derived category to a cosilting complex.
\ed{Pro}

\Pf. $\Rightarrow$ \ \ Since short exact sequences give triangles in the derive category, it is easy to see that every cotilting module is cosilting in the derived category by Theorem \ref{ATTA}.

\ \hskip 17pt $\Leftarrow$\ \ Note that there is a faithful embedding from Mod$R$ into $\mc{D}$. i.e.,
for any two modules $M,N\in\mr{Mod}R$, we have that
$\HD(M,N)\cong \mr{Hom}_R(M,N)$. Moreover, we have $\HD(M,N[i])\cong \mr{Ext}_R^i(M,N)$ for all $i>0$ and for any two modules $M,N$. So the condition (2) in the definition of cotilting modules is satisfied.

As to the condition (1) in the definition of cotilting modules,
since $T$ is isomorphic in the derived category to a cosilting complex and $T\in \mc{D}^{\geq}$, by Theorem \ref{ATTA}, we have
that $T\in \widecheck{\AdpD{Q}}$. i.e., there are triangles
$T_{i}\to^{\alpha_{i}} Q_{i}\to T_{i+1}\to$ with $Q_{i}\in \AdpD{Q}$ for all $0\leq i\leq n$, where $T_{0}=T$, $T_{n+1}=0$.
We will show that these triangles are in $\mr{Mod}R$ and hence give short exact sequences in $\mr{Mod}R$.

Consider firstly the triangle $T_{0}\to^{\alpha_{0}} Q_{0}\to T_{1}\to$, where $T_{0}=T$ is already an R-module. Then $\alpha_0\in\HD{(T_0,Q_0)}\simeq\mr{Hom}_R(T_0,Q_0)$ shows that $\alpha_0$ is homomorphism between modules. We claim that $\alpha_0$ is injective. To see this,
taking any momomorphism $\beta$: $T_0\to Q'$ with $Q'$ an injective module.
 Since $T\in \widecheck{\AdpD{Q}}$, it is easy to see that all $T_{i}\in$ $^{\bot_{i>0}}Q$.
Hence $\mathrm{Hom}_{R}(\alpha_{0},Q')\cong \HD(\alpha_{0},Q')$ is surjective.
Then we have the following commutative diagram in $\mr{Mod}R$ for some homomorphisms $\beta'$.

\vskip 10pt

 \setlength{\unitlength}{0.09in}
 \begin{picture}(50,8)

 \put(27,0){\makebox(0,0)[c]{$Q'$}}

 \put(12,6){\makebox(0,0)[c]{$0$}}
                              \put(13,6){\vector(1,0){2}}

 \put(18,6){\makebox(0,0)[c]{$ker(\alpha_{0})$}}
                              \put(21,6){\vector(1,0){3}}
                              \put(21,6){\vector(1,0){3}}

  \put(22,7){\makebox(0,0)[c]{$i$}}

 \put(26,3){\makebox(0,0)[c]{$\beta$}}

  \put(31,7){\makebox(0,0)[c]{$\alpha_{0}$}}

 \put(27,6){\makebox(0,0)[c]{$T_{0}$}}
                              \put(30,6){\vector(1,0){3}}
                              \put(27,4){\vector(0,-1){2}}
\put(34,3){\makebox(0,0)[c]{$\beta'$}}
                              \put(34,5){\vector(-1,-1){5}}

 \put(36,6){\makebox(0,0)[c]{$Q_{0}$}}
\end{picture}
\vskip 10pt

From the diagram above, we have that $\beta\alpha_0=\beta'\alpha_{0}i=0$. Note that $\beta$ is injective, so we obtain that $i=0$ and consequently, $\alpha_{0}$ is injective. Then we have an exact sequence $0\to T_{0}\to Q_{0}\to \mr{Coker}(\alpha_{0})\to 0$ which induces
a triangle $T_{0}\to Q_{0}\to \mr{Coker}(\alpha_{0})\to$. It follows that $T_{1}\cong \mr{Coker}(\alpha_{0})$
is (quasi-isomorphic to) an $R$-module. Repeating discussion above for all $i$, we can get that
each $\alpha_{i}$ is injective and each $T_{i}$ is (quasi-isomorphic to) an $R$-module.

Note that $T_{n}=Q_{n}$. By discussion above, we can get a long exact sequence
$0\to T\to Q_{0}\to Q_{1} \to \cdots \to Q_{n}\to0$. So $\mr{id}T\leq n$, i.e.,  the condition (1) in the definition of cotilting modules is satisfied.

Finally, still by Theorem \ref{ATTA}, we have that $Q\in \widehat{\AdpD{T}}$.
Similarly to the above process, we can get a long exact sequence
$0\to T_{n}\to \cdots \to T_{1} \to T_{0}\to Q\to 0 $ with $T_{i} \in \mathrm{Adp}T$, i.e.,  the condition (3) in the definition of cotilting modules is satisfied.
\ \hfill$\Box$


\vskip 15pt

\bg{Pro}\label{STTD}%
Suppose that $T\in \mathcal{D}^{\geq}$ is a cosilting complex.

$(1)$ If $S\bigoplus T$ is also a cosilting complex for some $S$, then
$S\in \AdpD{T}$.

\hangafter 1 \hangindent 35pt
$(2)$ If there are triangles $T_{i}\to Q_i\to T_{i+1}\to$ with $Q_i\in \AdpD{Q}$ for all $0\leq i\leq n$,
where $T_{0}=T$, $T_{n+1}=0$, then $\AdpD({\bigoplus_{i=0}^{n}Q_{i}})=\AdpD{Q}$.

$(3)$ If there are triangles $Q_{i+1}\to T_i\to Q_i\to$ with $T_i\in \AdpD{T}$ for all $0\leq i\leq m$,
where $Q_{0}=Q$, $Q_{m+1}=0$, then $\bigoplus_{i=0}^{m}T_{i}$ is a cosilting complex. Moreover,
$\AdpD({\bigoplus_{i=0}^{m} T_{i}})=\AdpD{T}$.
\ed{Pro}

\Pf. (1) Since $S\bigoplus T$ is a cosilting complex, we have that
$\langle \AdpD{(S\bigoplus T)}\rangle=\mc{K}^b(\mr{Inj}R)=\langle \AdpD{T}\rangle$.
It is easy to verify that $S\in$ $^{\bot_{i>0}}T$ and $S \in (\AdpD{T})^{\bot_{i>0}}$ since $S\bigoplus T$
is prod-semi-selforthogonal. It follows from Proposition \ref{OIIE} that

$\ \ \ \ \ \ \ \ \ \ \ \ \ \ \ \ \ \ \  \  $ $S\in$ $^{\bot_{i>0}}T\bigcap \langle\AdpD{(S\bigoplus T)}\rangle=$ $^{\bot_{i>0}}T\bigcap \langle \AdpD{T}\rangle=\widecheck{\AdpD{T}}$.
$\ \ \ \ \ \ \ \ \ \  \ \ \ \ \ \  \ \ \ \  $
Hence $S\in \widecheck{\AdpD{T}}\bigcap (\AdpD{T})^{\bot_{i>0}}=\AdpD{T}$ by Proposition \ref{TIRA}.

(2) Obviously, ${\bigoplus_{i=0}^{n}Q_{i}}\in \AdpD{Q}$ is prod-semi-selforthogonal.
It is not difficult to see that $\AdpD{T}\subseteq  \langle{\bigoplus_{i=0}^{n}Q_{i}}\rangle$ by Proposition \ref{AACU}. It follows that
$$\langle \AdpD{Q} \rangle=\mc{K}^b(\mr{Inj}R)=\langle \AdpD{T}\rangle
\subseteq \langle \AdpD({\bigoplus_{i=0}^{n}Q_{i}})\rangle\subseteq \langle \AdpD{Q} \rangle.$$
Hence $\langle \AdpD({\bigoplus_{i=0}^{n}Q_{i}})\rangle=\langle \AdpD{Q} \rangle$.
Clearly, $(\bigoplus_{i=0}^{n}Q_{i})\bigoplus Q$ is cosilting. It follows from (1) that
$\AdpD({\bigoplus_{i=0}^{n}Q_{i}})=\AdpD{Q}$.

(3) It is easy to verify that both $({\bigoplus_{i=0}^{m} T_{i}})\oplus T$ and ${\bigoplus_{i=0}^{m} T_{i}}$
are cosilting by Theorem \ref{ATTA}.
Consequently, we have that
$\AdpD({\bigoplus_{i=0}^{m} T_{i}})=\AdpD{T}$ by (1).
\ \hfill$\Box$

\bg{Pro}\label{cs-mn}%
Let $T\in \mathcal{D}^{\geq}$ be a cosilting complex and $n\geq 0$. Then
$T\in (\widecheck{\AdpD{Q}})_n$ if and only if $Q\in
(\widehat{\AdpD{T}})_n$.
\ed{Pro}
\Pf. $\Rightarrow$ We have that $Q\in (\widehat{\AdpD{T}})_m$ for some $m$, by Theorem \ref{ATTA}. If $m\leq n$, then the conclusion holds clearly. Suppose that $m> n$. There are triangles $Q_{i+1}\to T_{i}\to Q_{i}\to$ with $T_{i} \in \AdpD{T}$
for $0\leq i\leq m$, where $Q_{0}=Q$, $Q_{m+1}=0$. Applying the functor $\mathrm{Hom}_{\mathcal{D}}(-,Q_{m})$ to these triangles, we can obtain that
\bg{verse}

$\HD(Q_{m-1},Q_{m}[1]) \simeq \HD(Q_{m-2},Q_{m}[2])$

\ \hskip 101pt   $ \simeq \cdots \simeq\HD(Q_{0},Q_{m}[m])=\HD(Q,Q_{m}[m])$.

\ed{verse}
It is not difficult to verify that $\AdpD{T}\subseteq (^{\bot_{i>0}}Q)^{\bot_{i>n}}$ since
$T\in (\widecheck{\AdpD{Q}})_n$. Then we have that $\HD(Q,Q_{m}[t])=0$ for $t>n$ since $Q \in$ $^{\bot_{i>0}}Q$ and $Q_{m}=T_{m}\in \AdpD{T}$. Consequently, $\HD(Q_{m-1},Q_{m}[1])=0$ and  the triangle $T_{m}=Q_{m}\to T_{m-1}\to Q_{m-1}\to$ is split.
Hence $Q_{m-1}\in \AdpD{T}$ and $Q\in (\widehat{\AdpD{T}})_{m-1}$. By continuing this process,
we can finally obtain that $Q\in (\widehat{\AdpD{T}})_n$.

$\Leftarrow$ The proof is just the dual of above statement.
\ \hfill$\Box$

\vskip 10pt

The following is an easy observation.

\bg{Lem}\label{STTI}%

Suppose that $T\in\mc{D}$ is prod-semi-selforthogonal, then
$^{\bot_{i>0}}T={\mc{X}}_{\AdpD{T}}$.

\ed{Lem}

\Pf. Clearly, ${\mc{X}}_{\AdpD{T}}\subseteq$ $^{\bot_{i>0}}T$.

Take any $M\in$ $^{\bot_{i>0}}T$ and consider the triangle
$M\to^{\alpha} T^{X}\to M_{1}\to$, where $\alpha$ is the canonical evaluation map.
Applying the functor $\mathrm{Hom}_{\mathcal{D}}(-,T)$ to this triangle, we can obtain that
$\mathrm{Hom}_{\mathcal{D}}(M_{1},T[i])=0$ for all $i>0$. i.e., $M_{1}\in$ $^{\bot_{i>0}}T$. Continuing this process,
we get triangles  $M_{j}\to T_{j}\to M_{j+1}$ with $T_{j}\in \AdpD{T}$ and $M_{j}\in$ $^{\bot_{i>0}}T$
for all $j\geq0$, where $M_{0}=M$. Consequently, $M \in {\mc{X}}_{\AdpD{T}}$ by the
definition. So $^{\bot_{i>0}}T\subseteq{\mc{X}}_{\AdpD{T}}$
and the conclusion holds.
\ \hfill$\Box$

\bskip\

We say a complex is partial cosilting, if it satisfies
the first two conditions in Definition \ref{ACTI}.

\bg{Pro}\label{ITIP}%
If $T\in \mathcal{D}^{\geq}$ is partial cosilting. Then $T$ is cosilting
if and only if $^{\bot_{i>0}}T\subseteq \mathcal{D}^{\geq}$.
\ed{Pro}

\Pf. $\Rightarrow$ By Theorem \ref{ATTA} (3), there are triangles $Q_{i+1}\to T_{i}\to Q_{i}\to$ with
$T_{i}\in \AdpD{T}$ for all $0\leq i\leq n$, where $Q_{0}=Q, Q_{n+1}=0$.
Applying the functor $\mathrm{Hom}_{\mathcal{D}}(M,-)$ to these triangles, where $M\in$ $^{\bot_{i>0}}T$,
we can get that $M\in$ $^{\bot_{i>0}}Q$. Hence $^{\bot_{i>0}}T\subseteq$ $^{\bot_{i>0}}Q=\mathcal{D}^{\geq}$.

$\Leftarrow$ It is not difficult to verify that $Q \in$ $^{\bot_{i>n}}T$ for some $n$,
so $Q[-n]\in$ $^{\bot_{i>0}}T$.
Then there are triangles $Q[-i-1]\to 0\to Q[-i]\to$
for all $0\leq i\leq n-1$. Hence $Q \in (\widehat{^{\bot_{i>0}}T})_{n}$.
We can obtain a triangle $Q\to X\to Y\to $
with $X\in (\widehat{\AdpD{T}})_{n}, Y\in$ $^{\bot_{i>0}}T$ by Corollary \ref{FALT} and Lemma \ref{STTI}.
Note that $^{\bot_{i>0}}T\subseteq \mathcal{D}^{\geq}$. So $Y\in$ $^{\bot_{i>0}}Q$ and
the triangle $Q\to X\to Y\to$ is split. Hence $Q\in (\widehat{\AdpD{T}})_{n}$ since $(\widehat{\AdpD{T}})_{n}$
is closed under direct summands. Consequently, $T$ is a cosilting complex.
\ \hfill$\Box$

\vskip 10pt

We say that a complex $T\in \mathcal{D}^{\geq}$ is $n$-cosilting if it is
a cosilting complex such that $Q\in (\widehat{\AdpD{T}})_n$. A characterization of $n$-cotilting modules says that an $R$-module $T$ is $n$-cotilting if and only if $\mr{Cogen}^nT=\mr{KerExt}_R^{i>0}(-,T)$, see \cite{Bazz-nt} for more details. We will present a similar characterization of $n$-cosilting complexes in the following.

We need the following subcategory of $\mathcal{D}$.
 Let $T\in \mc{D}$ and $n>0$, we denote

\bg{verse}

$\mr{Copres}_{\mathcal{D}^{\geq}}^n(\AdpD{T})=\{M\in\mc{D}\ |$ there exist some triangles
$M_{i}\to T_i\to M_{i+1}\to$ with $T_i\in \AdpD{T}$ for all $0\le i<n$, where $M_n\in \mathcal{D}^{\geq}$
and $M_{0}=M$ $\}$.

\ed{verse}

It is not difficult to verify that $\mr{Copres}_{\mathcal{D}^{\geq}}^n(\AdpD{T})$ is closed under
products. The following result gives more properties about this subcategory.

\bg{Lem}\label{SAIR}%
$(1)$ $\mathcal{D}^{\geq}[-n]\subseteq \mr{Copres}_{\mathcal{D}^{\geq}}^n(\AdpD{T})$.

\ \hskip 45pt $(2)$ If $T\in \mathcal{D}^{\geq}$, then $\mr{Copres}_{\mathcal{D}^{\geq}}^n(\AdpD{T})\subseteq \mathcal{D}^{\geq}$.
\ed{Lem}

\Pf. Since $0\in \AdpD{T}$ and $\mathcal{D}^{\geq}$ is resolving, it is easy to verify the conclusions by the definitions.
\ \hfill$\Box$

\bg{Pro}\label{ITDI}%
Assume that $T\in \mathcal{D}^{\geq}$ is $n$-cosilting. Then $^{\bot_{i>0}}T=\mr{Copres}_{\mathcal{D}^{\geq}}^n(\AdpD{T})$.

\ed{Pro}

\Pf. By Proposition \ref{ITIP} and Lemma \ref{STTI}, we get that $^{\bot_{i>0}}T\subseteq$ $^{\bot_{i>0}}Q=\mathcal{D}^{\geq}$ and $^{\bot_{i>0}}T={\mc{X}}_{\AdpD{T}}$.
In particular, $^{\bot_{i>0}}T\subseteq\mr{Copres}_{\mathcal{D}^{\geq}}^n(\AdpD{T})$.

Now we prove that $\mr{Copres}_{\mathcal{D}^{\geq}}^n(\AdpD{T})\subseteq$ $^{\bot_{i>0}}T$.
For any $M\in \mathcal{D}^{\geq}=$ $^{\bot_{i>0}}Q$, it is not difficult to verify that $M\in$ $^{\bot_{i>n}}T$.
Take any $N\in \mr{Copres}_{\mathcal{D}^{\geq}}^n(\AdpD{T})$, then there are triangles
$N_{i}\to T_i\to N_{i+1}\to$ with $T_i\in \AdpD{T}$ for all $0\le i<n$, where $N_n\in \mathcal{D}^{\geq}$
and $N_{0}=N$. Applying the functor $\mathrm{Hom}_{\mathcal{D}}(T,-)$ to these triangles, we have that

\bg{verse}
$\mathrm{Hom}_{\mathcal{D}}(N_{0},T[i])\cong \mathrm{Hom}_{\mathcal{D}}(N_{1},T[i+1]) \cong \cdots \cong \mathrm{Hom}_{\mathcal{D}}(N_{n},T[i+n])=0$, $i>0.$

\ed{verse}
Hence $\mathrm{Hom}_{\mathcal{D}}(N,T[i])=\mathrm{Hom}_{\mathcal{D}}(N_{0},T[i])=0$, $i>0$ and $N\in$ $^{\bot_{i>0}}T$.
i.e., $\mr{Copres}_{\mathcal{D}^{\geq}}^n(\AdpD{T})\subseteq$ $^{\bot_{i>0}}T$.
The proof is then completed.
\ \hfill$\Box$
\mskip\

It is well known that $T\in \mc{K}^b(\mr{Inj}R)$ if and only if $\mathcal{D}^{+}\subseteq$
$^{\bot_{i\gg 0}}T$.

\bg{Pro}\label{ASSU}%
Assume that $T\in \mathcal{D}^{\geq}$. If $^{\bot_{i>0}}T=\mr{Copres}_{\mathcal{D}^{\geq}}^n(\AdpD{T})$,
then $T$ is $n$-cosilting.

\ed{Pro}

\Pf. Note that $T\in  \mr{Copres}_{\mathcal{D}^{\geq}}^n(\AdpD{T})=$$^{\bot_{i>0}}T$.
Since $\mr{Copres}_{\mathcal{D}^{\geq}}^n(\AdpD{T})$ is closed under
products, we have that $\AdpD{T}\subseteq$ $^{\bot_{i>0}}T$. Hence $T$ is prod-semi-selforthogonal.

From Lemma \ref{SAIR}, we know that $\mathcal{D}^{\geq}[-n]\subseteq \mr{Copres}_{\mathcal{D}^{\geq}}^n(\AdpD{T})
=$ $^{\bot_{i>0}}T$. So $\mathcal{D}^{\geq}\subseteq$ $^{\bot_{i>n}}T$. In particular, $\mathcal{D}^{+}\subseteq$
$^{\bot_{i\gg 0}}T$. Consequently, $T\in \mc{K}^b(\mr{Inj}R)$.
It follows that $T\in (\widecheck{\AdpD{Q}})_{m}$ for some $m$ from the argument above and Proposition \ref{OIIE}.

We note that $^{\bot_{i>0}}T=\mr{Copres}_{\mathcal{D}^{\geq}}^n(\AdpD{T})\subseteq \mathcal{D}^{\geq}$
by Lemma \ref{SAIR} (2). Hence, by Proposition \ref{ITIP}, $T$ is $m$-cosilting complex,  where  $m$ is the integer given in the last paragraph.
Since $T\in \widecheck{\AdpD{Q}}$. then there are triangles
$T_{i}\to Q_i\to T_{i+1}\to$ with $Q_i\in \AdpD{Q}$ for all $0\le i\leq m$, where $T_{m+1}=0$
and $T_{0}=T$. Applying the functor $\mathrm{Hom}(Q_{m},-)$ to these triangles, we can obtain that

\vskip 5pt

$\HD(Q_{m},T_{m-1}[1]) \cong \HD(Q_{m},T_{m-2}[2]) \cong \cdots \cong\HD(Q_{m},T_{0}[m])=\HD(Q_{m},T[m]).$

\vskip 5pt

Noted that $Q[-n]\in \mr{Copres}_{\mathcal{D}^{\geq}}^n(\AdpD{T})=$$^{\bot_{i>0}}T$
since $0\in \AdpD{T}$. Hence $\HD(Q,T[i+n])=0$, for any $i>0$.
If $m \leq n$, then $T$ is clearly $n$-cosilting. If $m>n$, it follows from the discussion  above
that $T_{m-1}\to Q_{m-1}\to Q_{m}\to$ is split. i.e., $T\in (\widecheck{\AdpD{Q}})_{m-1}$.
Repeating this process, we finally get that $T\in (\widecheck{\AdpD{Q}})_{n}$.
Consequently, $T$ is $n$-cosilting.
\ \hfill$\Box$
%

\vskip 10pt
Combining Proposition \ref{ITDI} and Proposition \ref{ASSU},
we obtain the following characterization of $n$-cosilting complexes.

\bg{Th}\label{c2}%
Assume that $T\in \mathcal{D}^{\geq}$. Then the following are equivalent:

$(1)$ $T$ is $n$-cosilting,

$(2)$ $^{\bot_{i>0}}T=\mr{Copres}_{\mathcal{D}^{\geq}}^n(\AdpD{T})$.
\ed{Th}

%
%
%


\hskip 15pt

%

In \cite{AR}, Auslander and Reiten
showed that there is a one-one correspondence between
isomorphism classes of basic cotilting  modules and certain contravariantly
finite  resolving  subcategories. Extending this result, Buan \cite{Buan} showed that there is a one-to-one correspondence between basic cotilting
complexes and certain contravariantly finite subcategories of the bounded derived category
of an artin algebra. In the following, we aim to extend such a result to
cosilting complexes.

We need the following definitions. Let $\mathcal{X}\subseteq\mathcal{Y}$ be two subcategories of $\mc{D}$.
$\mathcal{X}$ is said to be contravariantly finite in $\mathcal{Y}$,
if for any $Y\in \mathcal{Y}$, there is a homomorphism $f$ : $X\to Y$ for some $X\in \mathcal{X}$ such that
$\HD(X',f)$ is surjective for any $X'\in \mathcal{X}$. Moreover, $\mathcal{X}$ is said to be specially contravariantly finite in $\mathcal{Y}$,
if for any $Y\in \mathcal{Y}$, there is triangle $U\to X\to Y\to$ with some $X\in \mathcal{X}$ such that
$\HD(X',U[1])=0$ for any $X'\in\mathcal{X}$. Note that in the later case, one has that $U\in \mathcal{X}^{\bot_{i>0}}$ if $\mc{X}$ is closed under $[-1]$.

\bg{Pro}\label{ATTM}%
Assume that $T\in \mathcal{D}^{\geq}$ is cosilting. Then
$\widehat{^{\bot_{i>0}}T}=\mathcal{D}^{+}$ and $^{\bot_{i>0}}T\subseteq \mathcal{D}^{\geq}$  is specially contravariantly
finite in $\mathcal{D}^{+}$.
\ed{Pro}

\Pf. We have proved that $^{\bot_{i>0}}T\subseteq \mathcal{D}^{\geq}$ in Proposition \ref{ITIP},
so we get that $\widehat{^{\bot_{i>0}}T}\subseteq\mathcal{D}^{+}$.
Now we take any $X\in \mathcal{D}^{+}$. It is easy to see that $X\in$ $^{\bot_{i>m}}T$ for some $m$
since $T \in \mc{K}^b(\mr{Inj}R)$. Consequently, $X[-m]\in$ $^{\bot_{i>0}}T$.
Note that $0\in \AdpD{T}$, so we have $X\in \widehat{^{\bot_{i>0}}T}$ by the definition.
Hence $\widehat{^{\bot_{i>0}}T}=\mathcal{D}^{+}$.

By Lemma \ref{STTI}, we have that $^{\bot_{i>0}}T={\mc{X}}_{\AdpD{T}}$.
Taking any $X \in \mathcal{D}^{+}=\widehat{^{\bot_{i>0}}T}=\widehat{{\mc{X}}_{\AdpD{T}}}$,
by Corollary \ref{FALT}, we obtain a triangle $U\to V\to X\to$ with $U\in \widehat{\AdpD{T}}$ and $V\in$ $^{\bot_{i>0}}T$
. Note that $\widehat{\AdpD{T}}\subseteq (^{\bot_{i>0}}T)^{\bot_{i>0}}$, so particularly we get that  $\HD(M,U[1])=0$ for any $M\in$ $^{\bot_{i>0}}T$. It follows that
$^{\bot_{i>0}}T$ is specially contravariantly
finite in $\mathcal{D}^{+}$.
\ \hfill$\Box$
%

\mskip\

\bg{Pro}\label{ATIS}%
Assume that $\mc{T}\subseteq \mathcal{D}^{\geq}$ is specially contravariantly finite
in $\mathcal{D}^{+}$ and is resolving such that $\widehat{\mathcal{T}}=\mathcal{D}^{+}$. If
$\mc{T}\bigcap \mathcal{T}^{\bot_{i>0}}$ is closed under products,
then there is a cosilting complex $T$ such that
$\mc{T}={^{\bot_{i>0}}T}$.
\ed{Pro}

\Pf. It is not difficult to verify that
$\mathcal{D}^{+}=\widehat{\mathcal{T}}\subseteq$ $^{\bot_{i\gg0}}(\mathcal{T}^{\bot_{i>0}})$.
Hence we can obtain that $\mathcal{T}^{\bot_{i>0}}\subseteq \mc{K}^b(\mr{Inj}R)$.

Taking any $M\in \mathcal{D}^{+}$, since $\mc{T}$ is specially contravariantly finite
in $\mathcal{D}^{+}$ and is resolving, there are triangles
$M_{j+1}\to T_{j}\to M_{j}\to$ with $T_{j}\in \mathcal{T}$ for all  $j\geq0$, where $M_{0}:=M$
and each $M_{j}\in \mathcal{T}^{\bot_{i>0}}$ for $j\geq1$. It follows that $T_{j}\in \mc{T}\bigcap \mathcal{T}^{\bot_{i>0}}$
for all $j\geq1$.
Since $\mathcal{T}^{\bot_{i>0}}\subseteq \mc{K}^b(\mr{Inj}R)$ and $M\in$ $\mathcal{D}^{+}$,
it is easy to see that $M \in$ $^{\bot_{i>n}}(\mathcal{T}^{\bot_{i>0}})$ for some $n$ depending on $M$.
Applying $\mathrm{Hom}_{\mathcal{D}}(-,M_{n+1})$ to the triangles above, we obtain that $\mathrm{Hom}_{\mathcal{D}}(M_{n},M_{n+1}[1])\simeq \cdots\simeq \mathrm{Hom}_{\mathcal{D}}(M,M_{n+1}[n+1])=0$.
Thus, the triangle $M_{n+1}\to T_{n}\to M_{n}\to$ is split and so $M_{n}$ is a direct summand of $T_{n}$.
Note that the conditions $\mc{T}\bigcap \mathcal{T}^{\bot_{i>0}}$ is closed under products and
$\mc{T}$ is resolving imply that $\mc{T}\bigcap \mathcal{T}^{\bot_{i>0}}$ is closed under direct summands, so we have that $M_{n}\in \mc{T}\bigcap \mathcal{T}^{\bot_{i>0}}$.

Recall that $Q$ is an injective cogenerator in $\mr{Mod}R$. Note that $Q\in \mathcal{T}^{\bot_{i>0}}$ since $\mc{T}\subseteq \mathcal{D}^{\geq}=$ $^{\bot_{i>0}}Q$. Specially the object $M$ in the above to be $Q$, we obtain  triangles
$Q_{j+1}\to T'_{j}\to Q_{j}\to$ with  $Q_{j}\in \mathcal{T}^{\bot_{i>0}}$ and $T'_{j}\in \mc{T}\bigcap \mathcal{T}^{\bot_{i>0}}$ for all  $0\leq j\leq n$, where $Q_{0}=Q$ and $Q_{n+1}=0$.
Taking $T=\bigoplus_{j=0}^{n}T'_{j}$. We will show that $T$ is cosilting. It is easy to see that $T$ is precosilting since
$\mathcal{T}^{\bot_{i>0}}\subseteq \mc{K}^b(\mr{Inj}R)$ and
$\mc{T}\bigcap \mathcal{T}^{\bot_{i>0}}$ is closed under products. Moreover, the argument above shows that $Q\in \widehat{\AdpD{T}}$ too.  Hence $T$ is cosilting.

Now we need only prove that  $\mc{T}={^{\bot_{i>0}}T}$.
Obviously, we have that $\mc{T}\subseteq {^{\bot_{i>0}}T}$ since $T=\bigoplus_{j=0}^{n}T'_{j}$ and
$T'_{j}\in \mc{T}\bigcap \mathcal{T}^{\bot_{i>0}}$ for all  $0\leq j\leq n$.
Taking any $N\in {^{\bot_{i>0}}T}$.
Similar to the discussion above, there are triangle
$N_{j+1}\to T''_{j}\to N_{j}\to$ with  $N_{j}\in \mathcal{T}^{\bot_{i>0}}$,  $T''_0\in \mc{T}$ and $T''_{j}\in \mc{T}\bigcap \mathcal{T}^{\bot_{i>0}}$ for all  $1\leq j\leq m$, where $N_{0}=N$ and $N_{m+1}=0$.
Note that all objects in these triangles are in ${^{\bot_{i>0}}T}$.
For any $L\in \mc{T}\bigcap \mathcal{T}^{\bot_{i>0}}$, it is easy to verify that $T\bigoplus L$
is also cosilting, hence $L\in \AdpD{T}$ by Proposition \ref{STTD} (1). It follows that $\mc{T}\bigcap \mathcal{T}^{\bot_{i>0}}\subseteq \AdpD{T}$.
Now it is easy to see that $\mc{T}\bigcap \mathcal{T}^{\bot_{i>0}}= \AdpD{T}$.
So the above triangles imply that $N_{1}\in \widehat{\AdpD{T}}$ and consequently, $N_{1}\in {^{\bot_{i>0}}T} \bigcap\widehat{\AdpD{T}}=\AdpD{T}$.
Thus, the triangle $N_{1}\to T_{0}\to N\to$ is split. It follows that $N\in \mc{T}$ from the facts that $T_{0}, N_{1} \in \mc{T}$
and that $\mc{T}$ is resolving. So we obtain that $^{^{\bot_{i>0}}T}\subseteq \mc{T}$.
The proof is then completed.
\ \ \ \ \ \ \ \hfill$\Box$
%

%
\mskip\

By Propositions \ref{ATTM} and \ref{ATIS}, we obtain the following desired result. Here, we say two complexes $M$ and $N$ are equivalent if $\AdpD{M}=\AdpD{N}$.

\bg{Th}\label{TIAO}%
There is a one-one correspondence, given by $u: T \mapsto{^{\bot_{i>0}}T}$,
between equivalent class of cosilting complexes in $\mathcal{D}^{\geq}$
and subcategories $\mc{T}\subseteq \mathcal{D}^{\geq}$ which is specially
contravariantly finite in $\mathcal{D}^{+}$, resolving and closed under
products such that $\widehat{\mathcal{T}}=\mathcal{D}^{+}$.

\ed{Th}

\Pf. It follows from Propositions \ref{ATTM} and \ref{ATIS} that the correspondence is well-defined. Moreover, $u$ is surjective by Proposition \ref{ATIS}.
If both $T_{1}$ and $T_{2}$ are cosilting with ${^{\bot_{i>0}}T_1}={^{\bot_{i>0}}T_2}$,
it is easy to varify that $T_{1}\bigoplus T_{2}$ is also cosilting by the definition.
So we have that $\AdpD{T_{1}}=\AdpD{T_{2}}$ by Proposition \ref{STTD}, i.e., $T_1$ and $T_2$ are equivalent.
Hence, $u$ is bijective.
\ \hfill$\Box$

\vskip 30pt


\section{Quasi-cotilting modules and cosilting modules}
\sskip\
\

In this section, we introduce cosilting modules which is the dual of silting modules introduced in \cite{AMV}. We study their relationship with quasi-cotilting modules and provide some characterizations of cosilting modules. In particular, we obtain that all cosilting modules are pure-injective and cofinendo and that every presilting module has a Bongartz complement.

Let $R$ be a ring and $\mc{U}$ be a class of $R$-modules. Following \cite{Bazz-nt}, we denote by $\mr{Gen}^n\mc{U}$ the class of all modules $M$ such that there is an exact sequence $U_n\to\cdots\to U_1\to M\to 0$ with each $U_i\in \mc{U}$. Note that $\mr{Gen}^1\mc{U}$ and $\mr{Gen}^2\mc{U}$ are often denoted by $\mr{Gen}\mc{U}$ and $\mr{Pres}\mc{U}$ respectively. Dually, we denote by $\mr{Cogen}^n\mc{U}$ the class of all modules $M$ such that there is an exact sequence $0\to M\to U_1\to\cdots\to U_n$ with each $U_i\in \mc{U}$. Also we have $\mr{Cogen}\mc{U}={\mr{Cogen}}^1\mc{U}$ and $\mr{Copres}\mc{U}={\mr{Cogen}}^2\mc{U}$.

We simply denote $\mr{Gen}^n\mc{U}$ with $\mr{Gen}^nT$ in case that $\mc{U}=\mr{Add}T$ for some $R$-module $T$, where $\mr{Add}T$ denotes the class of modules which is a direct summand of some direct sums of copies of $T$. Dually, we simply denote $\mr{Cogen}^n\mc{U}$ with $\mr{Cogen}^nT$ in case that $\mc{U}=\mr{Adp}T$ for some $R$-module $T$, where $\mr{Adp}T$ denotes the class of modules which is a direct summand of some direct products of copies of $T$. We have similar simple notations $\mr{Gen}T$, $\mr{Pres}T$, $\mr{Cogen}T$, $\mr{Copres}T$.

Now let $T$ be an $R$-module. Recall that $T$ is an $n$-star module if $\mr{Gen}^nT=\mr{Gen}^{n+1}T$ and $\mr{Hom}_R(T,-)$ preserve the exactness of  exact sequences in $\mr{Gen}^nT$ \cite{Wnt}. Dually, one call that $T$ is an $n$-costar module if $\mr{Cogen}^nT=\mr{Cogen}^{n+1}T$ and $\mr{Hom}_R(-,T)$ preserve the exactness of  exact sequences in $\mr{Cogen}^nT$ \cite{Cnst}.

Let $\mc{U}$ be a class of $R$-modules. We denote by $\mr{KerExt}_R^1(\mc{U},-)$ the class of all $R$-modules $M$ such that $\mr{Ext}_R^1({U},M)=0$ for all $U\in\mc{U}$. We have similar notations such as $\mr{KerExt}_R^1(-,\mc{U})$.

Recall that an $R$-module $E$ is Ext-injective in $\mc{U}$ if $E\in \mc{U}\bigcap\mr{KerExt}_R^1(\mc{U},-)$, and dually, an $R$-module $E$ is Ext-projective in $\mc{U}$ if $E\in \mc{U}\bigcap\mr{KerExt}_R^1(-,\mc{U})$.

%

\vskip 5pt

We have the following definitions  \cite{AMV, CDT, ZW}.

\bg{Def}\label{}%
$(1)$ An $R$-module $M$ is called a quasi-tilting module, if it is  a 1-star module and
 is Ext-projective in $\mathrm{Gen} M$.

$(2)$ An $R$-module $M$ is called a quasi-cotilting module, if it is  a 1-costar module and
 is Ext-injective in $\mathrm{Cogen} M$.
\ed{Def}

We will say that two quasi-cotilting modules $M, N$ are equivalent if $\mr{Adp}M=\mr{Adp}N$.

Let $Q$ be an injective cogenerator of $\mathrm{Mod}R$. Following \cite{ATT}, an $R$-module $M$ is called $Q$-cofinendo
if there exist a cardinal $\gamma$ and a map $f$: $M^{\gamma}\to Q $ such that for any
cardinal $\alpha$, all maps $M^{\alpha}\to Q$ factor through $f$. An $R$-module $M$ is cofinendo if there is some injective cogenerator $Q$ of $\mathrm{Mod}R$
such that $M$ is $Q$-cofinendo.

A short exact sequence $0\to A\to B\to C\to 0$ is called pure exact if the induced sequence $0\to N\otimes_{R}A\to N\otimes_{R}B\to N\otimes_{R}C\to 0$ is still exact  for any right $R$-module $N$. An $R$-module $M$ is called pure injective
if $\mathrm{Hom}_R(-,M)$ preserves the exactness of all pure exact sequences. We note that a remarkable properties of cotilting modules is that they are pure-injective \cite {Bazz-pi, Stv}.

Let $\mathcal{U}$ be a  class of $R$-modules and $N$ be an $R$-module. Recall that a homomorphism $f: U\to N$  is
called a precover, or a right $\mathcal{U}$-approximation, of $N$
if  $U\in\mc{U}$ and $\mathrm{Hom}_{R}(U',f)$ is surjective for any $U'\in \mathcal{U}$.
A $\mathcal{U}$-precover  $f: U\to N$  of $N$ is called a $\mathcal{U}$-cover, or a minimal right $\mathcal{U}$-approximation, of $M$ if any $g: U\to U$ such that $f=fg$ must be an isomorphism.
A $\mathcal{U}$-precover  $f: U\to N$  of $N$ is called special if $\mr{Ker}f\in\mr{KerExt}_R^1(\mc{U},-)$.
A class $\mathcal{U}$ of $R$-modules is said to be a precover class, or contrvariantly finite, provided that every $R$-module
has a $\mathcal{U}$-precover. Cover classes and special precover classes are defined similarly.

Recall that a class $\mathcal{U}$ of $R$-modules  is torsion-free if $\mathcal{U}$ is closed under direct products,
submodules and extensions, see \cite{Dtor}.

We collect some import results on quasi-cotilting modules from \cite{ZW} in the following proposition.

\bg{Pro}\label{Qct}\label{main24}\label{char24}\label{IARM}
\ \ Let $M$ be an $R$-module and $Q$ be an injective cogenerator of $\mathrm{Mod}R$.

$(1)$ If $\mathrm{Cogen} M\subseteq \mr{KerExt}_R^1(-,M)$, then $(\mr{KerHom}_R(-,M), \mathrm{Cogen}M)$ is a torsion pair.

$(2)$ All quasi-cotilting modules are pure-injective and cofinendo;

\hangafter 1 \hangindent 35pt
$(3)$ $M$ is quasi-cotilting module if and only if $M$ is Ext-injective in $\mathrm{Cogen} M$ and there is an exact sequence
$0\to M_{1}\to M_{0}\to^{\alpha} Q$ with $M_{0},M_{1}\in \mathrm{Adp} M$ such that $\alpha$ is a $\mathrm{Cogen} M$-precover.

$(4)$ $M$ is 1-cotilting if and only if $M$ is quasi-cotilting and $Q\in \mathrm{Gen}(\mathrm{Cogen} M)$.

$(5)$ There are one-one correspondences between the following three classes

\begin{adjustwidth}{30pt}{0cm}

 \quad \ \ (i) equivalent classes of quasi-cotilting modules,

 (ii) torsionfree cover classes and,

 (iii) torsionfree specially precover classes.

\end{adjustwidth}

\ed{Pro}


We now turn to cosilting modules.

Firstly, let us recall the definition of silting modules given in \cite{AMV}. Let $\sigma: P_1\to P_0$ be in $\mr{Proj}R$ and $\mc{D}_{\sigma}$ be the class of all $R$-modules $N$ such that $\mr{Hom}_R(\sigma, N)$ is surjective. Then an $R$-module $M$ is said to be presilting if there is some $\sigma$ in $\mr{Proj}R$ such that $M=\mr{Ker}\sigma\in\mc{D}_{\sigma}$ and that $\mc{D}_{\sigma}$ is a torsion class. Moreover, an $R$-module $M$ is  said to be  silting if there is some $\sigma$ in $\mr{Proj}R$ such that $M=\mr{Ker}\sigma$ and $\mc{D}_{\sigma}=\mr{Gen}M$. Note that silting modules are always presilting \cite{AMV}.

As for the duall case, we take a homomorphism $\sigma: E_0\to E_1$ in $\mathrm{Inj}R$, and consider the associated class of $R$-modules

\hskip  110pt $\mathcal{F}_{\sigma}=\{M\in$ Mod$R$ $\mid$ $\mathrm{Hom}_{R}(M,\sigma)$ is surjective$\}$.

\vskip 10pt

The following result gives some useful properties of $\mathcal{F}_{\sigma}$. Here we say that a homomorphism $\sigma: E_0\to E_1$ in $\mathrm{Inj}R$ is an injective copresentation of $M$ if $M\simeq \mr{Ker}\sigma$.

\bg{Lem}\label{fsigma}%
Let $\sigma$ be a  homomorphism in $\mathrm{Inj}R$ with $K=\mr{Ker}\sigma$.
Then the following assertions hold:

$(1)$ $\mathcal{F}_{\sigma}$ is closed under submodules, extensions and direct sums;

$(2)$ $\mathcal{F}_{\sigma}\subseteq \mr{KerExt}_R^1(-,K)$;

\hangafter 1 \hangindent 35pt
$(3)$ An $R$-module $M$ belongs to $\mathcal{F}_{\sigma}$ if and only if $\mathrm{Hom}_{\mathcal{D}}(\theta[-1],\sigma)=0$, for any injective copresentation $\theta$ of $M$ (here $\theta, \sigma$ are considered as a complex with terms concentred on 0-th and 1-th positions).
\ed{Lem}

\Pf. (1) Easily.

(2) Factor $\sigma$: $E_{0}\to E_{1}$ canonically as $\sigma=i\pi$ with $i$: $C \to E_{1}$
and $\pi$: $E_{0}\to C$, where $C=\mathrm{Im} \sigma$.
For any $M\in \mathcal{F}_{\sigma}$, applying $\mathrm{Hom}_{R}(M,-)$ to the exact sequence $0\to K\to E_{0}\to^{\pi} C\to 0$
, we have an induced exact sequence
$\mathrm{Hom}_{R}(M,E_{0})\to^{\mathrm{Hom}_{R}(M,\pi)} \mathrm{Hom}_{R}(M,C)\to \mathrm{Ext}^{1}_{R}(M,K)\to 0$.
It is easy to verify that $\mathrm{Hom}_{R}(M,\pi)$ is surjective since that $\mathrm{Hom}_{R}(M,\sigma)$ is surjective implies $\mathrm{Hom}_{R}(M,i)$ is isomorphic.
So $\mathrm{Ext}^{1}_{R}(M,K)=0$.

(3) ($\Rightarrow$) Set $\theta$: $I_{0}\to I_{1}$ to be an injective copresentation of $M$ and write $\theta=i\pi$ with $i$: $\mathrm{Im} \theta \to I_{1}$
and $\pi$: $I_{0}\to \mathrm{Im} \theta$. For any morphism $f$: $I_{0}\to E_{1}$, We consider the following commutative diagram, where the second row is a complex.

 \setlength{\unitlength}{0.09in}
 \begin{picture}(50,10)

 \put(14,1){\makebox(0,0)[c]{$0$}}
                              \put(16,1){\vector(1,0){3}}
 \put(21,1){\makebox(0,0)[c]{$E_{0}$}}
                             \put(23,1){\vector(1,0){3}}
                             \put(25,2){\makebox(0,0)[c]{$_{\sigma}$}}

 \put(28,1){\makebox(0,0)[c]{$E_{1}$}}
                             \put(30,1){\vector(1,0){3}}

 \put(35,1){\makebox(0,0)[c]{$0$}}

                 \put(21,5.5){\vector(0,-1){3}}
             \put(22,4){\makebox(0,0)[c]{$_{g}$}}
                 \put(28,5.5){\vector(0,-1){3}}
                 \put(29,4){\makebox(0,0)[c]{$_{f}$}}

                 \put(26.5,4){\makebox(0,0)[c]{$_{s_{1}}$}}
                 \put(33.5,4){\makebox(0,0)[c]{$_{s_{0}}$}}

 \put(14,7){\makebox(0,0)[c]{$0$}}
                              \put(16,7){\vector(1,0){3}}

 \put(21,7){\makebox(0,0)[c]{$M$}}
                              \put(23,7){\vector(1,0){3}}
                              \put(25,8){\makebox(0,0)[c]{$_{\alpha}$}}

 \put(28,7){\makebox(0,0)[c]{$I_{0}$}}
                              \put(30,7){\vector(1,0){3}}
                              \put(32,8){\makebox(0,0)[c]{$_{\theta}$}}
                              \put(27,6){\vector(-1,-1){4}}
 \put(35,7){\makebox(0,0)[c]{$I_{1}$}}
                              \put(34,6){\vector(-1,-1){4}}

\end{picture}

 Since $M\in \mathcal{F}_{\sigma}$, we have a morphism $g$
such that $f\alpha=\sigma g$. There is a  morphism $s_{1}$ such that $g=s_{1}\alpha$ since
$E_{0}$ is injective.
And $(f-\sigma s_{1})\alpha=0$, thus, there is a morphism $h$: $\mathrm{Im} \theta \to E_{1}$ such that
$f-\sigma s_{1}=h\pi$. Since $E_{1}$ is injective and $i$ is monomorphic, we have a homomorphism $s_{0}$ such that $h=s_{0}i$.
It is easy to see that $f=\sigma s_{1}+s_{0}\theta$.

($\Leftarrow$) For any morphism $a$: $M\to E_{1}$, consider the following commutative diagram

 \setlength{\unitlength}{0.09in}
 \begin{picture}(50,10)

 \put(14,1){\makebox(0,0)[c]{$0$}}
                              \put(16,1){\vector(1,0){3}}
 \put(21,1){\makebox(0,0)[c]{$E_{0}$}}
                             \put(23,1){\vector(1,0){3}}
                             \put(25,1.5){\makebox(0,0)[c]{$_{\sigma}$}}

 \put(28,1){\makebox(0,0)[c]{$E_{1}$}}
                             \put(30,1){\vector(1,0){3}}

 \put(35,1){\makebox(0,0)[c]{$0$}}

                 \put(28,5.5){\vector(0,-1){3}}
                 \put(29,4){\makebox(0,0)[c]{$_{b}$}}

                 \put(26.5,4){\makebox(0,0)[c]{$_{s_{1}}$}}
                 \put(33.5,4){\makebox(0,0)[c]{$_{s_{0}}$}}

 \put(14,7){\makebox(0,0)[c]{$0$}}
                              \put(16,7){\vector(1,0){3}}

 \put(21,7){\makebox(0,0)[c]{$M$}}
                              \put(22,6){\vector(1,-1){4}}
                              \put(23,7){\vector(1,0){3}}
                              \put(25,8){\makebox(0,0)[c]{$_{\alpha}$}}
                              \put(24,5){\makebox(0,0)[c]{$_{a}$}}

 \put(28,7){\makebox(0,0)[c]{$I_{0}$}}
                              \put(30,7){\vector(1,0){3}}
                              \put(32,8){\makebox(0,0)[c]{$_{\theta}$}}
                              \put(27,6){\vector(-1,-1){4}}
 \put(35,7){\makebox(0,0)[c]{$I_{1}$}}
                              \put(34,6){\vector(-1,-1){4}}

\end{picture}

Since $E_{1}$ is injective, we have a morphism $b$ such that $a=b\alpha$. Hence we have $s_{1}$ and $s_{0}$
such that $b=\sigma s_{1}+s_{0}\theta$ by the assumption that $\mathrm{Hom}_{\mathcal{D}}(\theta[-1],\sigma)=0$.
It is not difficult to verify that $a=\sigma s_{1}\alpha$. Thus
$M\in \mathcal{F}_{\sigma}$.
\ \hfill $\Box$

\vskip 10pt

We also need the following result. The proof is simple, so we  left to the reader.

\bg{Lem}\label{Fab}%

$(1)$ Let $a$ and $b$ be a morphism in $\mathrm{Inj} R$, then $\mathcal{F}_{a\oplus b}=\mathcal{F}_{a}\bigcap \mathcal{F}_{b}$.

$(2)$ Let $\alpha$: $I_{0}\to I$ and $\beta$: $I_{1}\to I$ be two morphisms in $\mathrm{Inj} R$,
then $\mathcal{F}_{\alpha}\subseteq \mathcal{F}_{(\alpha,\beta)}$.

\ed{Lem}

We now give the definition of cosilting modules.

\bg{Def}\label{cosilting modules}%
$(1)$ An $R$-module $M$ is called precosilting if there is an injective copresentation $\sigma$ of M
such that $M\in \mathcal{F}_{\sigma}$ and $\mathcal{F}_{\sigma}$ is a torsion-free class.

$(2)$  An $R$-module $M$ is  called cosilting if there is an injective copresentation $\sigma$ of M
such that $\mathcal{F}_{\sigma}=\mathrm{Cogen} M$.
\ed{Def}


\bg{Rem}\label{remark}%
$(1)$ If $M$ is precosilting, then $\mathrm{Cogen} M\subseteq\mathcal{F}_{\sigma}\subseteq  \mr{KerExt}_R^1(-,M)$
by Lemma \ref{fsigma};

$(2)$ If $M$ is cosilting, then $\mathrm{Cogen} M=\mathcal{F}_{\sigma}\subseteq \mr{KerExt}_R^1(-,M)$. In particular, we have a torsion pair
$(\mr{KerHom}_R(-,M), \mathrm{Cogen}M)$ by Proposition \ref{IARM}. So
all cosilting modules are precosilting.
\ed{Rem}

We say  that an $R$-module $M$ is cosincere if $\mathrm{Hom}(M,I)\neq 0$ for any $0\neq I \in \mathrm{Inj}R$.

Recall that  an $R$-module $M$ is said to
1-cotilting, if it satisfies the following three conditions:
(i) $\mr{id}M \leq 1$, i.e., the injective dimension of $M$ is not more than 1,
(ii) $\mathrm{Ext}^{1}_{R}(M^{X},M)=0$ for every set $X$, and
(iii) there is an exact sequence $0\to M_1\to M_0\to Q\to
0$ with $M_0, M_1\in \mathrm{Adp}{M}$, where $Q$ is some injective cogenerator.
An $R$-module $M$ is called partial 1-cotilting if it just satisfies the first two conditions above.

\bg{Pro}\label{cotilting and cosilting}%
$(1)$ An R-module M is partial 1-cotilting (resp., 1-cotilting) if and only if M is a  precosilting (resp., cosilting)
module with respect to a surjective injective copresentation.

$(2)$ Suppose that $\mr{id}M \leq 1$. Then $M$ is 1-cotilting if and only if $M$ is a cosincere cosilting module.

\ed{Pro}

\Pf. (1) If $M$ is a 1-precotilting module, then $\mr{id}M \leq 1$. So we have a short exact sequence
$0 \to M\to I_{0}\to^{\sigma} I_{1}\to 0$ with $I_{0}$ and $I_{1}$ injective. It is easy to verify that $\mathcal{F}_{\sigma}= \mr{KerExt}_R^1(-,M)$ in the case.
Since $M$ is 1-precotilting, we have $M \in  \mr{KerExt}_R^1(-,M)=\mathcal{F}_{\sigma}$. So $M$ is precosilting.
If $M$ is 1-cotilting, we have $\mathrm{Cogen} M= \mr{KerExt}_R^1(-,M)=\mathcal{F}_{\sigma}$, thus $M$ is cosilting.
The converse is similar.

(2) It is easy to see that all cotilting modules are cosincere. By (1), all 1-cotilting modules are also cosilting.

Assume that $M$ is cosilting with respect to some $\sigma: I_0\to I_1$ in $\mr{Inj}R$, we have an exact sequence $0\to M\to I_{0}\to^{\sigma} I_{1}\to^{\pi} C\to 0$
with $I_{0}$ and $I_{1}$ injective.
Set $K=\mathrm{Im} \sigma$, then $K$ is injective since $\mr{id}M \leq 1$. It follows that the exact sequence
$0 \to K\to I_{1}\to C\to 0$ is split. Thus, $C$ is injective.
For any morphism $g$: $M\to C$, there is a morphism $f$: $M\to I_{1}$ such that $g=\pi f$.
Note that $f$ factors  through $\sigma$ since $M\in\mc{F}_{\sigma}$, we have  $g=0$. So $C=0$ since $M$ is cosincere.
Cosequently, $M$  is 1-cotilting by (1).
\ \hfill $\Box$

\vskip 10pt

The following result gives some relations among cosilting modules, quasi-cotilting modules and 1-cotilting modules. In particular, It shows all cosilting modules are pure-injective and cofinendo.

\bg{Pro}\label{acmaq}%
$(1)$ All cosilting modules are quasi-cotilting. In particular, all cosilting modules are pure-injective and cofinendo.

$(2)$ Let M be an R-module and Q be an injective cogenerator. If $Q\in \mathrm{Gen}(\mathrm{Adp} M)$,
then the following statements are equivelent:

(i) M is cotilting;

(ii) M is cosilting;

(iii) M is quasi-cotilting.

\ed{Pro}
\Pf. (1) Let $M$ be a cosilting module with respect to a homomorphism $\sigma$ in $\mr{Inj}R$. Then $\mathrm{Cogen} M=\mathcal{F}_{\sigma}\subseteq\mr{KerExt}_R^1(-,M)$.
Then we only need to prove that $\mathrm{Cogen}M\subseteq\mathrm{Copres}M$ by the definition.
For any $T \in \mathrm{Cogen}M$, we have a short exact sequence
$0\to T\to^{u} M^{X}\to C\to 0$ with $u$ the canonical evaluation map. It is enough to prove that
$C \in \mathcal{F}_{\sigma}=\mathrm{Cogen}M$. For any morphism $f$: $C\to I_{1}$,
we consider the following commutative diagram:

 \setlength{\unitlength}{0.09in}
 \begin{picture}(50,10)

 \put(14,1){\makebox(0,0)[c]{$0$}}
                              \put(16,1){\vector(1,0){3}}
 \put(21,1){\makebox(0,0)[c]{$M$}}
                             \put(23,1){\vector(1,0){3}}
                             \put(25,2){\makebox(0,0)[c]{$_{i}$}}

 \put(28,1){\makebox(0,0)[c]{$I_{0}$}}
                             \put(30,1){\vector(1,0){3}}
                             \put(32,2){\makebox(0,0)[c]{$_{\sigma}$}}

 \put(35,1){\makebox(0,0)[c]{$I_{1}$}}

                 \put(21,5.5){\vector(0,-1){3}}
             \put(22,4){\makebox(0,0)[c]{$_{h}$}}
                 \put(28,5.5){\vector(0,-1){3}}
                 \put(29,4){\makebox(0,0)[c]{$_{g}$}}
                 \put(35,5.5){\vector(0,-1){3}}
                 \put(36,4){\makebox(0,0)[c]{$_{f}$}}
                 \put(26.5,4){\makebox(0,0)[c]{$_{\alpha}$}}
                 \put(33.5,4){\makebox(0,0)[c]{$_{\beta}$}}

 \put(14,7){\makebox(0,0)[c]{$0$}}
                              \put(16,7){\vector(1,0){3}}

 \put(21,7){\makebox(0,0)[c]{$T$}}
                              \put(23,7){\vector(1,0){3}}
                               \put(25,8){\makebox(0,0)[c]{$_{u}$}}

 \put(28,7){\makebox(0,0)[c]{$M^{X}$}}
                              \put(30,7){\vector(1,0){3}}
                              \put(32,8){\makebox(0,0)[c]{$_{\pi}$}}
                              \put(27,6){\vector(-1,-1){4}}
 \put(35,7){\makebox(0,0)[c]{$C$}}
                              \put(36,7){\vector(1,0){3}}
                              \put(34,6){\vector(-1,-1){4}}

 \put(41,7){\makebox(0,0)[c]{$0$}}

\end{picture}

Since $M^{X} \in \mathcal{F}_{\sigma}$, there is a morphism
$g$ such that $f\pi=\sigma g$. As $u$ is the canonical evaluation map, we have a
morphism $\alpha$ such that $h=\alpha u$.
It is easy to verify that $(g-i\alpha)u=0$. Thus there exists $\beta$ such that
$g-i\alpha=\beta \pi$. And then $f\pi=\sigma g=\sigma \beta \pi$, $f=\sigma\beta$. Consequently,
$C \in \mathcal{F}_{\sigma}$. In particular, we get that all cosilting modules are pure-injective and cofinendo,
 by Proposition \ref{main24}.

(2) By Proposition \ref{cotilting and cosilting}, (1) and  Proposition \ref{main24}.
\ \hfill $\Box$

\vskip 10pt

We now give some characterizations of cosilting modules.

\bg{Pro}\label{precosilting and cosilting}%
Let $M$ be an $R$-module and $Q$ be an injective cogenerator. The following conditions are equivalent:

$(1)$  $M$ is a cosilting module with respect to $\sigma$.

\hangafter 1 \hangindent 35pt
$(2)$  $M$ is a precosilting module with respect to $\sigma$ and  there exists an exact sequence $0\to M_{1}\to M_{0}\to^{\varphi} Q$
with $M_{0}$ and $M_{1}$ in $\mathrm{Adp} M$ such that $\varphi$ is an $\mathcal{F}_{\sigma}$-precover.

$(3)$  $ \mathrm{Cogen} M\subseteq \mathcal{F}_{\sigma}$ and  there exists an exact sequence $0\to M_{1}\to M_{0}\to^{\varphi} Q$
with $M_{0}$ and $M_{1}$ in $\mathrm{Adp} M$ such that $\varphi$ is an $\mathcal{F}_{\sigma}$-precover.

\ed{Pro}

\Pf. $(3)\Rightarrow (1)$ Clearly, we only need to prove that $\mathcal{F}_{\sigma} \subseteq \mathrm{Cogen} M$. For any $T \in \mathcal{F}_{\sigma}$, there exists a monomorphism
$f$: $T\to Q^{Y}$ since $Q$ is an injective cogenerator. As $\varphi$ is an $\mathcal{F}_{\sigma}$-precover,
we have a morphism $g$: $T\to M_{0}^{Y}$ such that $f=\varphi^{Y}g$. Then $g$ is injective since $f$ is injective.
Thus $T \in \mathrm{Cogen} M$.

$(1)\Rightarrow (2)\Rightarrow (3)$ By Remark \ref{remark}, Proposition \ref{acmaq} and Proposition \ref{char24}.
\ \hfill $\Box$

\vskip 10pt

It is well known that  every partial 1-cotilting module
can be completed to a 1-cotilting module and the complement is usually called Bongartz complement.
The following result shows that every precosilting modules has also a  Bongartz complement.

\bg{Pro}\label{precosilting--cosilting}%

Every precosilting module $M$ with respect to an injective representation $\sigma$ is a direct summand of a cosilting
module $\widebar{M}=M\bigoplus N$ with same associated torsion-free class, that is, $\mathrm{Cogen} \widebar{M}=\mathcal{F}_{\sigma}$.

\ed{Pro}

\Pf. Set $\sigma$: $I_{0}\to I_{1}$. Taking the canonical evaluation map $u$: $Q\to I_{1}^{X}$ with $Q$ an injective cogenerator.
Consider the pullback diagram of $\sigma^{X}$ and $u$:

 \setlength{\unitlength}{0.09in}
 \begin{picture}(50,10)

 \put(14,1){\makebox(0,0)[c]{$0$}}
                              \put(16,1){\vector(1,0){3}}
 \put(21,1){\makebox(0,0)[c]{$M^{X}$}}
                             \put(23,1){\vector(1,0){3}}

 \put(28,1){\makebox(0,0)[c]{$I_{0}^{X}$}}
                             \put(30,1){\vector(1,0){3}}
                             \put(32,2){\makebox(0,0)[c]{$_{\sigma^{X}}$}}

 \put(35,1){\makebox(0,0)[c]{$I_{1}^{X}$}}

                 \put(21,5.5){\vector(0,-1){3}}

                 \put(28,5.5){\vector(0,-1){3}}
                 \put(27.5,4){\makebox(0,0)[c]{$_{v}$}}
                 \put(35,5.5){\vector(0,-1){3}}
                 \put(35.5,4){\makebox(0,0)[c]{$_{u}$}}

\put(14,7){\makebox(0,0)[c]{$0$}}
                              \put(16,7){\vector(1,0){3}}

 \put(21,7){\makebox(0,0)[c]{$M^{X}$}}
                              \put(23,7){\vector(1,0){3}}

 \put(28,7){\makebox(0,0)[c]{$N$}}
                              \put(30,7){\vector(1,0){3}}
                              \put(32,8){\makebox(0,0)[c]{$_{\phi}$}}

 \put(35,7){\makebox(0,0)[c]{$Q$}}

\end{picture}

Next, we prove that $N \in \mathcal{F}_{\sigma}$. For any morphism $f$: $N\to I_{1}$,
we consider the following commutative diagram:

 \setlength{\unitlength}{0.09in}
 \begin{picture}(50,10)

 \put(7,1){\makebox(0,0)[c]{$0$}}
                              \put(9,1){\vector(1,0){3}}
 \put(14,1){\makebox(0,0)[c]{$M$}}
                             \put(16,1){\vector(1,0){3}}

 \put(21,1){\makebox(0,0)[c]{$I_{0}$}}
                             \put(23,1){\vector(1,0){3}}
                             \put(25,2){\makebox(0,0)[c]{$_{\sigma}$}}

 \put(28,1){\makebox(0,0)[c]{$I_{1}$}}

                 \put(21,5.5){\vector(0,-1){3}}
             \put(21.5,4){\makebox(0,0)[c]{$_{g}$}}
                 \put(28,5.5){\vector(0,-1){3}}
                 \put(28.5,4){\makebox(0,0)[c]{$_{f}$}}

                 \put(26,4){\makebox(0,0)[c]{$_{h}$}}
                 \put(33,4){\makebox(0,0)[c]{$_{b}$}}

 \put(14,7){\makebox(0,0)[c]{$0$}}
                              \put(16,7){\vector(1,0){3}}

 \put(21,7){\makebox(0,0)[c]{$M^{X}$}}
                              \put(23,7){\vector(1,0){3}}
                               \put(25,7.5){\makebox(0,0)[c]{$_{\alpha}$}}

 \put(28,7){\makebox(0,0)[c]{$N$}}
                              \put(30,7){\vector(1,0){3}}
                              \put(32,7.5){\makebox(0,0)[c]{$_{\phi}$}}
                              \put(27,6){\vector(-1,-1){4}}
 \put(35,7){\makebox(0,0)[c]{$Q$}}

                              \put(34,6){\vector(-1,-1){4}}

\end{picture}

Similarly to discussion in the proof of Proposition \ref{acmaq} (1), we have that $f\alpha=\sigma g$, $g=h\alpha$ and
$f=\sigma h+b\phi$. From the pullback diagram above, we have $b\phi=\sigma v'$, where $b$ and $v'$ are component
maps $u$ and $v$ respectively. Thus $f=\sigma(h+v')$ and $N \in \mathcal{F}_{\sigma}$.

It is easy to see that $\phi$ is an $\mathcal{F}_{\sigma}$-precover from the universal property of the pullback.
Set $\widebar{M}=M\bigoplus N$. The pullback diagram above gives an exact sequence
$0\to N\to I_{0}^{X}\bigoplus Q\stackrel{(\sigma^X,u)}{\longrightarrow} I_{1}^{X}$. So we obtain an injective representation $\rho$ of $\widebar{M}$,
with $\rho=\sigma \oplus (\sigma^{X},u)$. Now we get that $\mathcal{F}_{\rho}=\mathcal{F}_{\sigma}$ and $\widebar{M}\in \mathcal{F}_{\rho}$, by Lemma \ref{Fab}. So $\widebar{M}$ is precosilting.
Combining with Proposition \ref{precosilting and cosilting}, we  know that $\widebar{M}$ is cosilting.
\ \hfill $\Box$


\bg{Rem}%
\emph{Most results in this section are also independently obtained by Breaz and Pop \cite{BP} recently}.
\ed{Rem}


\vskip 30pt

\section{AIR-cotilting modules}

\

In this section, we introduce AIR-cotilting modules and give precise relations between them and cosilting modules and quasi-cotilting modules. Moreover, it is shown that they are intimately related to 1-cosilting complexes. 

Let $R$ be a ring and $M$ be an $R$-module. Recall from \cite{Wlt} that $M$ is a {\bf large support $\tau$-tilting} module if it satisfies the following two conditions: (1) there is an exact sequence $P_1\to^{f} P_0\to T\to 0$ with $P_1, P_0$ projective such that $\mathrm{Hom}(f,T^{(X)})$ is surjective for any set $X$ and, (2) there is an exact sequence $R\to^g T_0\to T_1\to 0$ with $T_0, T_1\in \mathrm{Add}T$ such that $\mathrm{Hom}(g,T^{(X)})$ is surjective for any set $X$.

We will say that $M$ is an {\bf AIR-tilting} module if it is large support $\tau$-tilting.

It is easy to see that 1-tilting modules, support $\tau$-tilting modules over artin algebras \cite{AIR} and silting modules \cite{AMV} are all AIR-tilting modules. From the proof of the main theorem in \cite{Wlt}, we also know that an AIR-tilting module $M$ can always be completed to an equivalent silting module $\widebar{M}$ in sense that there is  some $M'\in\mr{Add}M$ such that $\widebar{M}=M\oplus M'$ is a silting module. It is known that both silting modules and AIR-tilting modules coincide with support $\tau$-tilting modules in the scope of the category of finitely generated modules over artin algebras. But it is a question if silting modules and AIR-tilting modules coincide with each other in general. It is also known that AIR-tilting modules are finendo quasi-tilting. But the converse is not true in general \cite{Trf}.


We introduce the following dual definition.

\bg{Def}\label{AIR-cotilting}%
 Let R be a ring and M be an R-module. M is called AIR-cotilting module if
 it satisfies the following conditions:

 \hangafter 1 \hangindent 35pt
$(1)$ there exists an exact sequence $0 \to M\to Q_{0}\to^{f} Q_{1}$
 such that $\mathrm{Hom}_{R}(M^{X},f)$ is surjective for any set $X$, where $Q_{0}$ and $Q_{1}$ are
 in $\mathrm{Inj}R$.

 $(2)$ there exists an exact sequence $0 \to M_{1}\to M_{0}\to^{g} Q$
 such that $\mathrm{Hom}_{R}(M^{X},g)$ is surjective for any set $X$, where Q is an injective
 cogenerator of $\mathrm{Mod}R$ and $M_{1}$, $M_{0} \in \mathrm{Adp} M$.

\ed{Def}%


An $R$-module $M$ will be called \textbf{partial AIR-cotilting} if it satisfies the first condition in the above definition.

\bg{Lem}\label{sublem}%

Let $\sigma: L\to Q_1$ be a homomorphism with $Q_1$ injective. Assume that $M$ is an $R$-module such that $\mr{Hom}_R(M,\sigma)$ is surjective. Then every submodule $N$ of $M$ has the property that $\mr{Hom}_R(N,\sigma)$ is surjective.

\ed{Lem}%


\Pf. Let $N$ be a submodule of $M$ and $\phi: N\to M$ be the canonical embedding. Take any $f\in\mr{Hom}_R(N,Q_1)$ and consider the following diagram.

 \setlength{\unitlength}{0.09in}
 \begin{picture}(50,10)

 \put(23,1){\makebox(0,0)[c]{$L$}}
                             \put(25,1){\vector(1,0){3}}
                             \put(26,1.5){\makebox(0,0)[c]{$_{\sigma}$}}

 \put(29,1){\makebox(0,0)[c]{$ Q_{1}$}}

                 \put(28,5.5){\vector(0,-1){3}}
                 \put(27.5,5){\makebox(0,0)[c]{$_{f}$}}

                 \put(33,4){\makebox(0,0)[c]{$_{b}$}}

 \put(28,7){\makebox(0,0)[c]{$N$}}
                              \put(30,7){\vector(1,0){3}}
                              \put(32,7.5){\makebox(0,0)[c]{$_{\phi}$}}
\put(35,7){\makebox(0,0)[c]{$M$}}

                              \put(34,6){\vector(-1,-1){4}}
                              \put(32,6){\vector(-2,-1){7}}
                                   \put(30,5.5){\makebox(0,0)[c]{$_{h}$}}

\end{picture}

Since $Q_1$ is injective, $f$ lifts to a homomorphism $b: M \to Q_1$ such that $f=b\phi$. By the assumption, $b$ further lifts to a homomorphism $h: M\to L$ such that $b=\alpha h$. Then $f=b\phi=\sigma h\phi$, i.e., $f$ factors through $\sigma$. Hence we see that  $\mathrm{Hom}_{R}(N,\sigma)$ is surjective.
\hfill $\Box$

\vskip 10pt

The following proposition gives a characterization of partial AIR-cotilting modules in terms of $\mc{F}_{\sigma}$.

\bg{Pro}\label{STTE}%

An $R$-module $M$ is partial AIR-cotilting  if and only if there is an exact sequence $0 \to M\to Q_{0}\to^{\sigma} Q_{1}$ such that $\mathrm{Cogen} M\subseteq \mc{F}_{\sigma}$. In particular, $\mathrm{Cogen} M\subseteq  \mathrm{KerExt}_{R}^{1}(-,M)$ if  $M$ is partial AIR-cotilting.

\ed{Pro}%

\Pf. $(\Rightarrow)$ Suppose that $M$ is partial AIR-cotilting, i.e., there exists an exact sequence $0 \to M\to Q_{0}\to^{\sigma} Q_{1}$
 such that $\mathrm{Hom}(M^{X},\sigma)$ is surjective for any set $X$, where $Q_{0}$ and $Q_{1}$ are
in $\mathrm{Inj}R$.   Let $N\in\mr{Cogen}M$. Then there a monomorphism $\phi: N\to M^{X}$ for some $X$. Now by Lemma \ref{sublem} we obtain that  $\mathrm{Hom}_{R}(N,\sigma)$ is surjective. Thus, $\mathrm{Cogen} M\subseteq \mc{F}_{\sigma}$.

$(\Leftarrow)$ Obviously.
\hfill $\Box$

\vskip 10pt

We have the following easy corollary, which implies that, for an $R$-module $M$ of injective dimension not more than 1, $M$ is partial AIR-cotilting if and only if $M$ is partial 1-cotilting.

\bg{Cor}\label{LMBA}%

Let M be an R-module and $\mr{id}M\leq 1$. Then the following statements are equivalent:

$(1)$ $M$ is partial AIR-cotilting;

$(2)$ $\mathrm{Cogen} M\subseteq \ker \mathrm{Ext}_{R}^{1}(-,M)$;

$(3)$ $\mathrm{Ext}_{R}^{1}(M^{X},M)=0$ for any set X.

\ed{Cor}%

\Pf. (1)$\Rightarrow$(2) By Proposition \ref{STTE}.

(2)$\Rightarrow$(3) Obviously.

(3)$\Rightarrow$(1) Since $\mr{id}M\leq 1$, there is a short exact sequence
$0 \to M\to Q_{0}\to^{\alpha} Q_{1}\to 0$ with $Q_{0}$ and $Q_{1}$
in $\mathrm{Inj}R$. Applying $\mathrm{Hom}_{R}(M^{X},-)$ to this short exact sequence, we have
$$0 \to \mathrm{Hom}_{R}(M^{X},M)\to \mathrm{Hom}_{R}(M^{X},Q_{0})\to \mathrm{Hom}_{R}(M^{X},Q_{1})\to \mathrm{Ext}_{R}^{1}(M^{X},M)=0.$$
So $\mathrm{Hom}_{R}(M^{X},\alpha)$ is surjective.
\ \hfill $\Box$


\vskip 10pt We also have the following result.

\bg{Pro}\label{}%

Let Q be an injective cogenerator of $\mathrm{Mod}R$ and M be an R-module such that $Q\in \mathrm{Gen} M$.

$(1)$ If  $M$ is partial AIR-cotilting, then $\mr{id}M\leq 1$.

$(2)$ If M is AIR-cotilting, then M is 1-cotilting.

\ed{Pro}%

\Pf. Note that there is an exact sequence
$M^{(Y)}\to^{\gamma} Q \to 0$, since $Q\in\mathrm{Gen} M$.

(1) By the definition, we have an exact sequence $0 \to M\to Q_{0}\to^{\alpha} Q_{1}$
 such that $\mathrm{Hom}_{R}(M^{X},\alpha)$ is surjective for any set $X$, where $Q_{0}$ and $Q_{1}$ are
 in $\mathrm{Inj}R$.
As $Q$ is an injective cogenerator,
$Q_{1}$ is a summand of $Q^{X}$ for some $X$, and then we have a canonical projective
$\pi$: $Q^{X} \to Q_{1}$.
 Hence we have a surjection $f=\pi \gamma^{X}$: $(M^{(Y)})^{X}\to Q_{1}$.
Consider the following commutative diagram,
where $i$ is  a canonical embedding.

 \setlength{\unitlength}{0.09in}
 \begin{picture}(48,14)

 \put(36,12){\makebox(0,0)[c]{$(M^{(Y)})^{X}$}}
                            \put(36,11){\vector(0,-1){7}}
                            \put(39,12){\vector(1,0){3}}
                            \put(35,11){\vector(-1,-1){7}}
 \put(37,9){\makebox(0,0)[c]{$_{f}$}}
  \put(43,7){\makebox(0,0)[c]{$_{g}$}}
   \put(39,8){\makebox(0,0)[c]{$_{h}$}}
   \put(33,8){\makebox(0,0)[c]{$_{hi}$}}

\put(41,13){\makebox(0,0)[c]{$_{i}$}}

 \put(46,12){\makebox(0,0)[c]{$(M^{Y})^{X}$}}

                            \put(46,11){\vector(-1,-1){8}}
                            \put(43,11){\vector(-2,-1){15}}

 \put(17,2){\makebox(0,0)[c]{$0$}}
                              \put(18,2){\vector(1,0){3}}

 \put(22,2){\makebox(0,0)[c]{$M$}}
                              \put(23,2){\vector(1,0){3}}

 \put(28,2){\makebox(0,0)[c]{$Q_{0}$}}
                              \put(30,2){\vector(1,0){4}}

 \put(32,3){\makebox(0,0)[c]{$_{\alpha}$}}

 \put(36,2){\makebox(0,0)[c]{$Q_{1}$}}

\end{picture}

 There is a morphism $g$ such that $f=gi$ since
$Q_{1}$ is injective. Following from the property of the morphism $\alpha$, there exists a morphism $h$ such that $g=\alpha h$.
Hence, $f=gi=(\alpha h)i$. Then we obtain that $\alpha$ is surjective since $f$ surjective.
Thus, $\mr{id}M\leq 1$.

(2) Since $M$ is  AIR-cotilting,  there exists an exact sequence $0 \to M_{1}\to M_{0}\to^{\beta} Q$
 such that $\mathrm{Hom}_{R}(M^{X},\beta)$ is surjective for any $X$, where $M_{1}$ and $M_{0}$ are in $\mathrm{Adp} M$.
Clearly, $\mathrm{Hom}_{R}(M^{(Y)},\beta)$ is also surjective in the case, so we have a morphism $\delta$: $M^{(Y)}\to M_{0}$ such that $\gamma=\beta\delta$. Thus $\beta$ surjective
since  $\gamma$ is surjective. It follows that
$M$ is 1-cotilting from (1), Proposition \ref{LMBA} and the definition of 1-cotilting modules.
\ \hfill $\Box$

\vskip 10pt

Next we will consider the relations between 2-term cosilting complexes (i.e., 1-cosilting complexes in Section 2) and AIR-cotilting modules. We need some preparations.

\bg{Lem}\label{leti}%

Let $I^{\bullet}$: $0\to I_{0}\to^{\alpha} I_{1}\to 0$ be a 2-term complex of injective
modules and
$J^{\bullet}$: $0\to J_{0}\to J_{1}\to \cdots\ $  be a complex. If $K=H^{0}(J^{\bullet})$,
then the following statements are equivalent:

$(1)$ $K\in\mc{F}_{\alpha}$, i.e., $\mathrm{Hom}_{R}(K,\alpha)$ is surjective;

$(2)$ $\mathrm{Hom}_{\mathcal{D}}(J^{\bullet},I^{\bullet}[1])=0$.

In particular, an $R$-module $K\in\mc{F}_{\alpha}$ if and only if $\mathrm{Hom}_{\mathcal{D}}(K,I^{\bullet}[1])=0$.

\ed{Lem}%

\Pf. The proof is dual to Lemma 3.4 in \cite{AIR}.
\ \hfill $\Box$

$ $

By Lemma \ref{leti}, we can easily obtain the following corollary:

\bg{Cor}\label{STTIA}%
Suppose that there is an exact sequence $0 \to M\to I_{0}\to^{\alpha} I_{1}$ with $I_{0}$ and $I_{1}$ in $\mathrm{Inj}R$.
Then the following statements are equivalent:

$(1)$ $\mathrm{Hom}_{R}(M^{X},\alpha)$ is surjective for any set $X$;

$(2)$ The complex $I^{\bullet}$: $0\to I_{0}\to I_{1}\to 0$ is partial cosilting.

\ed{Cor}%


The following result shows that we can obtain AIR-cotilting modules and cosilting modules from 2-term cosilting complexes.

\bg{Pro}\label{AIR-cosilting1}%

Let  the 2-term complex $I^{\bullet}$: $0\to I_{0}\to^{\alpha} I_{1}\to 0$ be cosilting and set $K= H^{0}(I^{\bullet})$. Then

$(1)$ $K$ is AIR-cotilting.

$(2)$ $K$ is a cosilting module.

\ed{Pro}%

\Pf. (1) By Corollary \ref{STTIA}, we only need to check the condition (2) in Definition \ref{AIR-cotilting}.
Since $I^{\bullet}$ is cosilting, there is a triangle $I^{\bullet}_{1}\to I^{\bullet}_{0}\to^{\beta^{\bullet}} Q\to$
with $I^{\bullet}_{1}$ and $I^{\bullet}_{0}$ in $\mathrm{Adp}I^{\bullet}$, where $Q$ injective
cogenerator (see Section 2, Theorem \ref{ATTA} and Proposition \ref{cs-mn}).
By taking homologies, we can obtain an exact sequence
$0\to H^{0}(I^{\bullet}_{1})\to H^{0}(I^{\bullet}_{0})\to H^{0}(Q)(=Q)$.
Set $K_{1}=H^{0}(I^{\bullet}_{1})$, $K_{0}=H^{0}(I^{\bullet}_{0})$ and $\beta = H^{0}(\beta^{\bullet})$, then we have an exact sequence $0\to K_1\to K_0\to^{\beta} Q$. It is easy to see that $K_{1}$  and $K_{0}$ are in $\mathrm{Adp} K$.

For any $\gamma \in \mathrm{Hom}_{R}(K^{X},Q)$ with $X$ a set, we see that $\gamma$ lifts to a homomorphism $\gamma^{\bullet} \in \mathrm{Hom}_{\mathcal{D}}((I^{\bullet})^{X},Q)$,
since $Q$ is injective. By the assumption,  $I^{\bullet}$ is prod-semi-selforthogonal, so we have that $\HD{((I^{\bullet})^{X},\beta^{\bullet})}$ is surjective. Thus, there is a morphism $\delta^{\bullet}: (I^{\bullet})^{X}\to I^{\bullet}_{0}$
such that $\gamma^{\bullet}=\beta^{\bullet}\delta^{\bullet}$. Then we obtain that $H^{0}(\gamma^{\bullet})=H^{0}(\beta^{\bullet})H^{0}(\delta^{\bullet})$.
That is, $\gamma=\beta\delta$, where $\delta=H^{0}(\delta^{\bullet})$.
Hence $\mathrm{Hom}_{R}(K^{X},\beta)$ is surjective.

(2) It is easy to see that  $\mathrm{Cogen} K\subseteq  \mc{F}_{\alpha}$ by Corollary \ref{STTIA} and Lemma \ref{fsigma}. We claim that the exact sequence  $0\to K_1\to K_0\to^{\beta} Q$ obtained in (1) satisfies that $\beta$ is an $\mc{F}_{\alpha}$-precover. Thus $K$ is cosilting by Proposition \ref{precosilting and cosilting}.

In fact, take any $M\in \mc{F}_{\alpha}$ and any homomorphism $f: M\to Q$. We can consider these objects and homomorphisms in the derived category. By Lemma \ref{leti}, one has that $\HD{(M, I^{\bullet}[1])}=0$. Thus, by applying the functor $\HD{(M,-)}$ to the triangle $I^{\bullet}_{1}\to I^{\bullet}_{0}\to^{\beta^{\bullet}} Q\to$, we get that $\HD{(M, \beta^{\bullet})}$ is surjective, i.e., there exists some $\eta^{\bullet}: M\to I^{\bullet}_{0}$ such that $f=\beta^{\bullet}\eta^{\bullet}$.  Then we obtain that $H^{0}(f)=H^{0}(\beta^{\bullet})H^{0}(\eta^{\bullet})$.
That is, $f=\beta\eta$, where $\eta=H^{0}(\delta^{\bullet})$.
Hence $\mathrm{Hom}_{R}(M,\beta)$ is surjective and  $\beta$ is an $\mc{F}_{\alpha}$-precover.
\ \hfill $\Box$

\vskip 10pt

We will say that two AIR-cotilting modules $M, N$ are equivalent, denoted by $M\sim N$, provided that $\mr{Adp}M=\mr{Adp}N$.
The next result shows that  AIR-cotilting modules also give 2-term cosilting complexes.

\bg{Pro}\label{AIR-cosilting2}%
Let T be an AIR-cotilting module. Then

\hangafter 1 \hangindent 35pt
$(1)$ there is a 2-term cosilting complex $M^{\bullet}$ such that  $H^{0}(M^{\bullet})\simeq T\oplus T'$ with $T'\in\mr{Adp}T$. In particular, $H^{0}(M^{\bullet})\sim T$;

$(2)$ the cosilting complexs in $(1)$ is unique up to equivalences.

\ed{Pro}%

\Pf. (1) Since $T$ is an AIR-cotilting module, there exist two exact sequences
$0 \to T\to^{i} I_{0}\to^{\alpha} I_{1}$ and $0 \to T_{1}\to^{s} T_{0}\to^{t} Q$
such that $\mathrm{Hom}_{R}(T^{X},\alpha)$ and $\mathrm{Hom}_{R}(T^{X},t)$  are surjective respectively
for any set $X$, where $I_{0}$, $I_{1} \in \mathrm{Inj}R$ and $T_{1}$, $T_{0} \in \mathrm{Adp} T$.
Assume that $T_{0}'\oplus T_{0}=T^{Y}$ for some $R$-module $T_{0}'$,
then we have an exact sequence $0\to T_{1}'\to^{s'} T^{Y}\to^{t'} Q$
with $T_{1}'=T_{0}'\oplus T_{1}$, 
$s'=\jz{s}{0}{0}{1}$
and $t'=(t,0)$. It is easy to see that $\mathrm{Hom}_{R}(T^{X},t')$ is surjective.
Since $Q$ is injective, there exists a morphism $u$: $I_{0}^{Y}\to Q$ such that $t'=u\cdot i^{Y}$.
Set $I^{\bullet}$: $0\to I_{0}\to^{\alpha} I_{1}\to 0$. Clearly, $u$ induces a map of complex
$u^{\bullet}$: $(I^{\bullet})^{Y}\to Q$ with $H^{0}(u^{\bullet})=t'$.
Then we get a triangle $(I^{\bullet})^{Y}\to^{u^{\bullet}} Q\to \mathrm{Con}(u^{\bullet})\to$,
where $\mathrm{Con}(u^{\bullet})$: $0\to I_{0}^{Y}\to^{\theta} I_{1}^{Y}\bigoplus Q\to 0$ with
$\theta=\ljz{\alpha^Y}{u}$
is a complex with terms fixed in $(-1)$-th and $0$-th positions.
By taking homology of this triangle, we have an exact sequence
$0\to H^{-1}(\mathrm{Con}(u^{\bullet}))\to T^{Y}\to^{t'} Q$, thus
$H^{-1}(\mathrm{Con}(u^{\bullet}))\cong T_{1}'$.

We assert that $M^{\bullet}=I^{\bullet}\oplus \mathrm{Con}(u^{\bullet})[-1]$ is the desired cosilting complex.
Note that $H^{0}(M^{\bullet})=H^{0}(I^{\bullet})\bigoplus H^{-1}(\mathrm{Con}(u^{\bullet}))=T\bigoplus T_{1}'\sim T$,
since $T_{1}'\in \mathrm{Adp} T$. Obviously, $M^{\bullet}\in K^{b}(\mathrm{Inj} R)$ and
$Q\in \langle \mathrm{Adp}{M^{\bullet}}\rangle$. It remains to prove that $M^{\bullet}$  is prod-semi-selforthogonal
by Definition \ref{ACTI}. This is proceeded as follows.

(1) $\mathrm{Hom}_{\mathcal{D}}((M^{\bullet})^{X},I^{\bullet}[1])=0$ for any $X$. This is by Lemma \ref{leti}.

(2) $\mathrm{Hom}_{\mathcal{D}}((M^{\bullet})^{X},\mathrm{Con}(u^{\bullet})[-1][1])=0$ for any $X$.
Indeed, we only need to prove that $\mathrm{Hom}_{R}(K,\theta)$ is surjective by Lemma \ref{leti},
where $K=H^{0}((M^{\bullet})^{X})=(T\oplus T')^{X}$. Note that $K$ is a direct summand of $T^{X'}$ for some $X'$, it is sufficient to prove that  $\mathrm{Hom}_{R}(T^{X'},\theta)$ is surjective. For any morphism
$\ljz{f}{g}$:
$T^{X'}\to I_{1}^{Y}\bigoplus Q$.
Consider the following diagram, where $t', i^Y, \alpha^Y, u$ was defined as above. 

 \setlength{\unitlength}{0.09in}
 \begin{picture}(48,15)

 \put(17,11){\makebox(0,0)[c]{$0$}}
                              \put(17.5,11){\vector(1,0){3}}

 \put(22,11){\makebox(0,0)[c]{$T^{Y}$}}
                              \put(23,11){\vector(1,0){4}}
                              \put(21,10){\vector(0,-1){5}}

  \put(25,12){\makebox(0,0)[c]{$_{i^{Y}}$}}

 \put(29,11){\makebox(0,0)[c]{$I_{0}^{Y}$}}
                              \put(30,11){\vector(1,0){4}}
                              \put(28,10){\vector(-1,-1){5}}

 \put(32,12){\makebox(0,0)[c]{$_{\alpha^{Y}}$}}

 \put(36,11){\makebox(0,0)[c]{$I_{1}^{Y}$}}

 \put(37,4){\makebox(0,0)[c]{$T^{X'}$}}
                        \put(36,5){\vector(0,1){4}}
                        \put(35,4.8){\vector(-1,1){5}}
                        \put(35,4.2){\vector(-2,1){12}}
                        \put(35,4){\vector(-2,0){12}}

 \put(21,4){\makebox(0,0)[c]{$Q$}}

  \put(29,3){\makebox(0,0)[c]{$_{g}$}}
  \put(20,7){\makebox(0,0)[c]{$_{t'}$}}
  \put(24,7){\makebox(0,0)[c]{$_{u}$}}
  \put(29,6.5){\makebox(0,0)[c]{$_{b}$}}
  \put(37,7){\makebox(0,0)[c]{$_{f}$}}
  \put(34,7){\makebox(0,0)[c]{$_{a}$}}

\end{picture}

Since $\mathrm{Hom}_{R}(T^{X'},\alpha^{Y})$ is surjective, we can obtain a morphism $a$ such that
$f=\alpha^{Y}a$. 
Further, since $\mathrm{Hom}_{R}(T^{X'},t')$ is surjective, we can obtain some morphism $b, c: T^{X'}\to T^Y$ such that
$g=t'b$ and  $ua=t'c$. 
Note that $t'=u\cdot i^{Y}$ and $\alpha^Y i^Y=0$, it is easy to see that $\ljz{f}{g}=\ljz{\alpha^Y}{u}(a-i^{Y}c+i^{Y}b)$.

Thus $M^{\bullet}$ is just the desired cosilting complex.

(2) Suppose that $M^{\bullet}$ and $N^{\bullet}$ are 2-term cosilting complexes satisfying $H^{0}(M^{\bullet})\sim H^{0}(N^{\bullet})$. Then $\mathrm{Hom}_{\mathcal{D}}((M^{\bullet})^{X},N^{\bullet}[1])=0$
and $\mathrm{Hom}_{\mathcal{D}}((N^{\bullet})^{X},M^{\bullet}[1])=0$,
by Lemma \ref{leti}. Now it is easy to verify that $M^{\bullet}\bigoplus N^{\bullet}$ is a
cosilting complex, therefore, we have that $M^{\bullet} \in \AdpD N^{\bullet}$ and $N^{\bullet} \in \AdpD M^{\bullet}$ by Proposition \ref{STTD}. It follows that $\AdpD N^{\bullet}=\AdpD N^{\bullet}$, i.e., $M^{\bullet}$ and $N^{\bullet}$ are equivalent.
\ \hfill $\Box$

\vskip 10pt

As a direct corollary, we obtain the following relation between cosilting modules and AIR-cotilting modules, dual to the tilting case.

\bg{Cor}\label{cos-AIR}%

Let $M$ be AIR-cotilting. Then there is some $M'\in \mr{Adp}M$ such that $\widebar{M}=M\oplus M'$ is a cosilting module.

\ed{Cor}%


\bg{Lem}\label{itcu}%
If two cosilting complexes $U^{\bullet}$ and $V^{\bullet}$ are equivalent, then
$\mathrm{Adp}H^{k}(U^{\bullet})=\mathrm{Adp}H^{k}(V^{\bullet})$ for any integer k.

\ed{Lem}%

\Pf. It is enough  to show that $\mathrm{Adp}H^{k}(U^{\bullet})\subseteq\mathrm{Adp}H^{k}(V^{\bullet})$.
For any $M \in \mathrm{Adp}H^{k}(U^{\bullet})$, we have $M\bigoplus N=[H^{k}(U^{\bullet})]^{X}
\cong H^{k}((U^{\bullet})^{X})$ for some $N$. Since $\AdpD U^{\bullet}=\AdpD V^{\bullet}$,
we have $(U^{\bullet})^{X}\bigoplus W^{\bullet}=(V^{\bullet})^{Y}$ for some $W^{\bullet}$.
Thus, $H^{k}((U^{\bullet})^{X})\bigoplus H^{k}(W^{\bullet})\cong H^{k}((V^{\bullet})^{Y})$
and $M\bigoplus N \bigoplus H^{k}(W^{\bullet})\cong[H^{k}(V^{\bullet})]^{Y}$.
So  $M \in \mathrm{Adp}H^{k}(V^{\bullet})$.
\ \hfill $\Box$

\mskip\
Now we obtain the following theorem.

\bg{Th}\label{AIRcs}%
Let $R$ be a ring. There is a bijective correspondence between the equivalent classes
of AIR-cotilting modules and two-terms cosilting complexes.

\ed{Th}%

\Pf. By Proposition \ref{AIR-cosilting1}, \ref{AIR-cosilting2}
and Lemma \ref{itcu}.
\ \hfill $\Box$

\vskip 10pt

We now turn to study the precise relations between quasi-cotilting modules, cosilting modules and AIR-cotilting modules. To this aim, the following is a key result.

\bg{Lem}\label{cogenm}%
An $R$-module $M$ is partial AIR-cotilting if and only if $\mathrm{Cogen}M\subseteq \mr{KerExt}_R^1(-,M)$.

\ed{Lem}%

\Pf. ($\Rightarrow$) By Proposition \ref{STTE}.
%

($\Leftarrow$) Let $0\to M\to E_{0}\to^{f} E_{1}\to E_{2}$ be the minimal injective resolution of $M$, where each $E_i$ is in $\mathrm{Inj}R$. We need only to prove that $\mathrm{Hom}_{R}(M^{X},f)$ is surjective for any $X$.

Take any morphism $g$: $M^{X}\to E_{1}$ and consider
the following commutative diagram.

 \setlength{\unitlength}{0.09in}
 \begin{picture}(50,10)

 \put(7,1){\makebox(0,0)[c]{$0$}}
                                \put(9,1){\vector(1,0){3}}
 \put(14,1){\makebox(0,0)[c]{$M$}}
                              \put(16,1){\vector(1,0){3}}
 \put(21,1){\makebox(0,0)[c]{$E_{0}$}}
                             \put(23,1){\vector(1,0){3}}
                             \put(25,2){\makebox(0,0)[c]{$_{f}$}}

 \put(28,1){\makebox(0,0)[c]{$E_{1}$}}
                             \put(30,1){\vector(1,0){3}}
                             \put(32,2){\makebox(0,0)[c]{$_{\theta}$}}

 \put(35,1){\makebox(0,0)[c]{$E_{2}$}}

                 \put(21,5.5){\vector(0,-1){3}}
             \put(22,4){\makebox(0,0)[c]{$_{a}$}}
                 \put(28,5.5){\vector(0,-1){3}}
                 \put(29,4){\makebox(0,0)[c]{$_{g}$}}
                 \put(35,5.5){\vector(0,-1){3}}
                 \put(36,4){\makebox(0,0)[c]{$_{1}$}}
                 \put(26.5,4){\makebox(0,0)[c]{$_{t}$}}
                 \put(33.5,4){\makebox(0,0)[c]{$_{n}$}}

 \put(14,7){\makebox(0,0)[c]{$0$}}
                              \put(16,7){\vector(1,0){3}}

 \put(21,7){\makebox(0,0)[c]{$K$}}
                              \put(23,7){\vector(1,0){3}}
                               \put(25,8){\makebox(0,0)[c]{$_{k}$}}

 \put(28,7){\makebox(0,0)[c]{$M^{X}$}}
                              \put(30,7){\vector(1,0){3}}
                              \put(32,8){\makebox(0,0)[c]{$_{\theta g}$}}
                              \put(27,6){\vector(-1,-1){4}}
 \put(35,7){\makebox(0,0)[c]{$E_{2}$}}

                              \put(34,6){\vector(-1,-1){4}}

\end{picture}

Set $K=\mr{Ker}(\theta g)$ and factor $f=h\pi$ canonically, where $h: \mathrm{Im} f\to E_{1}$ and $\pi: E_{0}\to \mathrm{Im} f$.
There exists a morphism $\alpha: K\to \mathrm{Im} f$ such that $gk=h \alpha$, since $\theta gk=0$.
As $K\in \mathrm{Cogen M} \subseteq  \mr{KerExt}_R^1(-,M)$, we have that $\mr{Hom}_R(K,\pi)$ is surjective. Hence, there is some $a: K\to E_0$ such that $\alpha=\pi a$. Then we get $gk=h \alpha=h\pi a=fa$.

Since $E_0, E_1$ are injective, a canonical argument shows that there are two morphisms $t$ and $n$ such that
$a=tk$ and $g=n\theta g+ft$. Setting $\beta=g-ft$, then we have that
$n\theta\beta=n\theta(g-ft)=n\theta g=g-ft=\beta$,
that is, $\beta=n\theta\beta$. Now we claim that $\mathrm{Im}h\bigcap \mathrm{Im}\beta=0$. Indeed, for any $e \in \mathrm{Im}h\bigcap \mathrm{Im}\beta$, we have $h(x)=e=\beta(y)=n\theta\beta(y)=n\theta h(x)=0$, since $\theta h=0$. Thus $e=0$.
Since $h$ is an injective envelope by assumption, we have that $\mathrm{Im}h$ is an essential submodule of $E_{1}$. This implies that $\mathrm{Im}\beta=0$, i.e.,  $\beta=0$. Then we have $g=ft$.
So $\mathrm{Hom}(M^{X},f)$ is surjective for any set $X$.
\ \hfill $\Box$

\vskip 10pt

The following is a direct corollary.

\bg{Cor}\label{pairc-ds}%

A direct summand of a partial AIR-cotilting module is again partial AIR-cotilting.

\ed{Cor}


The following result is well-known.

\bg{Lem}\label{suppose} \label{direct summand}%

$(1)$ Suppose that $f$: $M \to E$ is an injective envelope of M and $g$: $M \to E'$ is a monomorphism
with $E'$ injective. Then $E'\simeq E\oplus E''$  and $g\simeq \ljz{f}{0}$.

$(2)$ Let $0\to M\to^{\alpha} E_{0}\to^{\beta} E_{1}$ be a minimal injective resolution of M and $0\to M\to^{\delta} I_{0}\to^{\sigma} I_{1}$
be any injective resolution of M. Then $I_0\simeq E_{0}\oplus E_0'$  and $I_1\simeq E_1\oplus E_0'\oplus E_1'$ and, moreover, the complex $0\to M\to^{\delta} I_{0}\to^{\sigma} I_{1}$ is isomorphic to the direct sums of three complexes $0\to M\to^{\alpha} E_{0}\to^{\beta} E_{1}$, $0\to 0\to E_0'\to^{1_{E_0'}} E_0'$ and $0\to 0\to 0\to E_1'$.

\ed{Lem}%

The following result gives a characterization of partial AIR-cotilting modules in term of its minimal injective copresentation.

\bg{Pro}\label{minima resolution}%

Let $M$ be an $R$-module and $0\to M\to E_{0}\to^{\delta} E_{1}$ be its minimal injective copresentation. Then $M$ is partial AIR-cotilting if and only if $\mathrm{Cogen}M\subseteq  \mc{F}_{\delta}$.

\ed{Pro}%

\Pf. ($\Leftarrow$) By the definition.

($\Rightarrow$) If $M$ is partial AIR-cotilting, then  there exists an exact sequence $0\to M\to I_{0}\to^{\sigma} I_{1}$ such that
$\mathrm{Hom}_{R}(M^{X},\sigma)$ is surjective for any set X, where $I_{0}$ and $I_{0}$ are in $\mathrm{Inj}R$.  By Lemma \ref{direct summand}, we see that $\sigma: I_0\to I_1$ is isomorphic to the direct sums of $\beta: E_0\to E_1$ and $1_{E_0'}: E_0'\to E_0'$ and $0\to E_2'$ for some injective modules $E_0'$ and $E_1'$. Thus, $\delta$ is a direct summand of $\sigma$. It follows that $\mathrm{Hom}_{R}(M^{X},\delta)$ is surjective from the assumption that $\mathrm{Hom}_{R}(M^{X},\sigma)$ is surjective. Thus, $M^X\in\mc{F}_{\delta}$ for any $X$. It follows that $\mathrm{Cogen}M\subseteq  \mc{F}_{\delta}$ by Lemma \ref{fsigma}.
\ \hfill $\Box$

\vskip 10pt

Moreover, we have the following easy observation by Lemma \ref{direct summand} and the involved definition.

\bg{Lem}\label{Fsgmlem}%

Let $\sigma: I_0\to I_1$ be a homomorphism of injective $R$-modules. Assume that $\alpha: I_0\to I$ is the composition of the canonical map $\pi: I_0\to \mr{Im}\sigma$ and the injective envelope map $i: \mr{Im}\sigma\to I$. Then $\sigma$ is isomorphic to the direct sum of $\alpha$ and the zero map $0\to I'$ for some injective module $I'$. In the case, it holds that $\mc{F}_{\sigma}=\mc{F}_{\alpha}\bigcap \mr{KerHom}_R(-,I')$.

\ed{Lem}

Now we are in the position to give our main result in the paper.

\bg{Th}\label{31}%

Let $M$ be an $R$-module. Then the following statements are equivalent.

$(1)$ $M$ is AIR-cotilting.

$(2)$ $M$  is quasi-cotilting.

$(3)$ $M$ is   cosilting.

\ed{Th}%

\Pf. $(1) \Rightarrow (2)$ By Lemma \ref{cogenm}, we have that $\mathrm{Cogen} M\subseteq \mr{KerExt}_R^1(-,M)$. Then it is easy to see that $M$ is Ext-injective in $\mr{Cogen}M$. Since $M$ is AIR-cotilting,  we also have an exact sequence $0 \to M_{1}\to M_{0}\to^{g} Q$
 such that $\mathrm{Hom}_{R}(M^{X},g)$ is surjective, where $Q$ is an injective
 cogenerator and $M_{1}, M_{0}\in \mathrm{Adp} M$.
 It is easy to see that $g$ is a $\mathrm{Cogen} M$-precover by Lemma \ref{sublem}. Hence, $M$ is quasi-cotilting  by Proposition \ref{char24}.

$(2) \Rightarrow (1)$ By Lemma \ref{cogenm} and Proposition \ref{char24}.

$(3) \Rightarrow (2)$ By Proposition \ref{acmaq}.

$(1) \Rightarrow (3)$   By Corollary \ref{cos-AIR},
there exists some $M'\in\mr{Adp}M$ such that $\widebar{M}\simeq M\oplus M'$  is a cosilting module.
Let $\sigma, \alpha, \beta$ be the minimal injective copresentations of $\widebar{M}, M, M'$ respectively. Then we have that $\sigma=\alpha\oplus\beta$.
Now we will show that $\mathcal{F}_{\sigma}=\mathcal{F}_{\alpha}$.
Since $ M^{'} \in \mathrm{Adp} M$, there is a morphism $\gamma$ in $\mr{Inj}R$ satisfiing
$\beta\oplus\gamma=\alpha^{X}$ for some $X$, so
$\mathcal{F}_{\alpha}=\mathcal{F}_{\alpha^{X}}=\mathcal{F}_{\beta\oplus\gamma}=\mathcal{F}_{\beta}
\bigcap\mathcal{F}_{\gamma}$. In particular, $\mathcal{F}_{\alpha}\subseteq \mathcal{F}_{\beta}$.
Thus, $\mathcal{F}_{\sigma}=\mathcal{F}_{\alpha\oplus\beta}=\mathcal{F}_{\alpha}
\bigcap\mathcal{F}_{\beta}=\mathcal{F}_{\alpha}$.

Since $\widebar{M}$  is a cosilting module,  there is an injective copresentation $\eta$ of $\widebar{M}$ $: I_0\to I_1$, where $I_0, I_1$ are injective, such that $\mr{Cogen}\widebar{M}=\mc{F}_{\eta}$.  By Lemma \ref{Fsgmlem}, we have that $\eta\simeq \sigma\oplus (0\to I')$ for some injective $R$-module $I'$, and then $\mc{F}_{\eta}=\mc{F}_{\sigma}\bigcap \mr{KerHom}_R(-,I')$. Let $\xi=\alpha\oplus (0\to I')$, then $\xi$ is an injective copresentation of $M$. Note that $\mc{F}_{\eta}=\mc{F}_{\xi}$ by the argument above and that $\mr{Cogen} M=\mathrm{Cogen} \widebar{M}$, so we obtain that $\mathrm{Cogen} M=\mr{Cogen}\widebar{M}=\mc{F}_{\eta}=\mathcal{F}_{\xi}$,
i.e., $M$ is a cosilting module.
\ \hfill $\Box$


\vskip 10pt
Combining results in the Proposition \ref{Qct} and Theorems \ref{AIRcs} and \ref{31}, we obtain the following result.

\bg{Th}\label{m11th}%

There are bijections between

$(1)$ equivalent classes of AIR-cotilting (resp., cosilting, quasi-cotilting) modules,

$(2)$ equivalent classes of 2-term cosilting complexes,

$(3)$ torsion-free cover classes and,

$(4)$ torsion-free special precover classes.

\ed{Th}
{\small

}


\begin{thebibliography}{17}

\bibitem{AIR}{} T. Adachi, O. Iyama and I. Reiten, $\tau$-tilting theory, Compos. Math. 150 (3) (2014), 415-452.

\bibitem{ATIO}{} T. Aihara and O. Iyama, Silting mutation in triangulated categories, J. London Math. Soc.  85  (3) (2012), 633-668.

\bibitem{AC}{}  L. Angeleri-H\"{u}gel and F.U. Coelho, Infinitely generated tilting modules of finite projective dimension, Forum Math. 13
(2001), 239-250.

\bibitem{AMV}{} L. Angeleri-H\"{u}gel, F. Marks and J. Vit\'{o}ria, Silting modules, Int. Math. Res. Notices (2015), in press. (See also arxiv: 1405.2531v1).


\bibitem{ATT}{} L. Angeleri-H\"{u}gel, A. Tonolo and J. Trlifaj, Tilting preenvelopes and cotilting precovers, Alg. Rep. Theory, 4 (2) (2001), 155-170.

\bibitem{AR}{} M. Auslander and I. Reiten, Applications of contravariantly finite subcategories,
                             Adv. Math. 86 (1991), 111-152.


\bibitem{Bazz-nt}{} S. Bazzoni, A characterization of n-cotilting and n-tilting modules, J. Algebra
                             273 (2004), 359-372.

\bibitem{Bazz-pi}{} S. Bazzoni, Cotilting modules are pure injective, Proc.  Amer. Math. Soc. 131(12) (2003), 3665-3672.

\bibitem{BP}{} S. Breaz and F. Pop, Cosilting modules, arXiv: 1510.05098v1.


\bibitem{SBMB}{} S. Brenner and M. Butler, Generalizations of the Bernstein-Gelfand-Ponomarev reflection functors, Lecture Notes in Math. 832 (1980), 103-169 .

\bibitem{Buan}{} A.B. Buan, Subcategories of the derived category and cotilting complexes, Coll. Math. 88(1) (2001), 1-11.

\bibitem{BZ}{} A.B. Buan and Y. Zhou, A silting theorem, arXiv: 1503.06129v1.


\bibitem{CDT}{} R. Colpi, G. D’este, and A. Tonolo, Quasi-tilting modules and counter equivalences,
J. Algebra 191 (1997), 461-494.

\bibitem{RCJT}{} R. Colpi and J. Trlifaj, Tilting modules and tilting torsion theories, J. Algebra 178 (1995), 614-634.

\bibitem{Dtor}{} S. E. Dickson, A torsion theory for Abelian categories, Trans. Amer. Math. Soc. 121 (1) (1966), 223-235.

\bibitem{Hp-b}{} D. Happel, Triangulated Categories in the Representation Theory of Finite Dimension Algebras,
London Math. Soc. Lect. Note Ser. 119 (1988).

\bibitem{HR}{} D. Happel and C. M. Ringel, Tilted algebras, Tran. Amer. Math. Soc. 274 (1982), 399-443.
%
%

\bibitem{Cnst} D. He, $n$-Costar modules and $n$-cotilting modules, Master thesis (in Chinese) (2009), Northwest Normal University.


\bibitem{BKDV}{} B. Keller and D. Vossieck, Aisles in derived categories, Bull. Soc. Math. Belg. Sér. A 40 (2) (1988), 239-253.


\bibitem{KY}{} S. Koenig and D. Yang, Silting objects, simple-minded collections, t-structures and co-t-structures for finite-dimensional algebras,  Doc. Math. 19 (2014), 403-438.

\bibitem{MSSS}{} O. Mendoza Hern\'{a}ndez,  E. C. S\'{a}enz Valadez, V. Santiago Vargas and M. J. Souto Salorio, Auslander-Buchweitz context and
co-t-structures, Appl. Categ. Structures 21 (5) (2013), 417-440.


\bibitem{YMIY}{} Y. Miyashita, Tilting modules of finite projective dimension, Math. Z. 193 (1986), 113-146.

\bibitem{Ric}{}  M. Rickard, Morita theory for derived categories, J. Lon. Math. Soc. 39 (2) (1989), 436-456. %

\bibitem{Stv}{}  J. Stovick, All $n$-cotilting modules are pure-injective, Proc.  Amer. Math. Soc. 134 (7) (2006), 1891-1897.

\bibitem{Trf}{} J. Trlifaj, On $\tau$-rigid modules, manuscript 2014. %

\bibitem{Wnt} J. Wei, $n$-Star modules and $n$-tilting modules, J. Algebra 283 (2005), 711-722.

\bibitem{Wst} J. Wei, Semi-tilting complexes, Israel J. Mathematics 194 (2013), 871-893.

\bibitem{Wlt} J. Wei, Large support $\tau$-tilting modules, preprint (2014).

\bibitem{Wns}{} J. Wei, Z. Huang, W Tong. and J. Huang, Tilting modules of finite projective dimension
and a generalization of star modules, Journal of Algebra 268 (2003), 404-418.

\bibitem{ZW} P. Zhang and J. Wei, Quasi-cotilting modules and torsion-free classes, preprint (2015).






\end{thebibliography}
\end{document}